\documentclass[11pt,reqno]{amsart}%
\usepackage{amsmath, amsfonts, amssymb, amsthm, amscd, amsbsy}
\usepackage[usenames,dvipsnames,svgnames,x11names,hyperref]{xcolor}
\usepackage{geometry}
\usepackage{graphicx}
\usepackage{hyperref}
\hypersetup{backref=true,
pagebackref=true,
hyperindex=true,
colorlinks=true,
breaklinks=true,
urlcolor=NavyBlue,
linkcolor=Fuchsia,
bookmarks=true,
bookmarksopen=false,
filecolor=black,
citecolor=ForestGreen,
linkbordercolor=red
}
\newtheorem{theorem}{Theorem}[section]
\theoremstyle{plain}
\newtheorem{corollary}{Corollary}[section]
\newtheorem{lemma}{Lemma}[section]
\newtheorem{proposition}{Proposition}[section]
\newtheorem{remark}{Remark}[section]
\numberwithin{equation}{section}
\allowdisplaybreaks
\AtBeginDocument{\hypersetup{pdfborder={0 0 0.01}}}
\DeclareMathOperator{\arsinh}{arsinh}

\begin{document}
\title[Optimal AB-program for  Poincar\'e--Sobolev and logarithmic Sobolev inequalities]{ Optimal AB-program for Poincar\'e--Sobolev and logarithmic Sobolev inequalities on hyperbolic spaces}
\author{Anh Xuan Do}
\address{Anh Xuan Do: Department of Mathematics, University of Connecticut, Storrs, CT
06269, USA}
\email{anh.do@uconn.edu}
\author{Debdip Ganguly}
\address{Debdip Ganguly: Theoretical Statistics and Mathematics Unit\\
Indian Statistical Institute, Delhi Centre\\
S.J. Sansanwal Marg, New Delhi, Delhi 110016, India}
\email{debdip@isid.ac.in}
\author{Nguyen Lam}
\address{Nguyen Lam: School of Science and the Environment, Grenfell Campus, Memorial
University of Newfoundland, Corner Brook, NL A2H5G4, Canada}
\email{nlam@mun.ca}
\author{Guozhen Lu}
\address{Guozhen Lu: Department of Mathematics, University of Connecticut, Storrs, CT
06269, USA}
\email{guozhen.lu@uconn.edu}

\begin{abstract}
Though the sharp Poincar\'e inequality and the sharp Sobolev inequality on the hyperbolic space $\mathbb{H}^N$ are well known, whether one can subtract the full spectral gap $\left(\frac{N-1}{p}\right)^{p}$ of the $p$-Laplacian and keep the sharp Sobolev constant $S(N,p)$ of $\mathbb{R}^N$ at the same time has been known only in a few cases. For $p=2$, this was proved by Benguria, Frank and Loss in dimension three, and it was shown by Hebey to be false in dimensions $N\geq4$. For $p\neq2$, it was proved by Nguyen for $N\geq4$ and $\frac{2N}{N-1}\leq p<N$. By employing, among other tools, a suitable change of function that maps the radial $p$-Dirichlet integral on $\mathbb{H}^N$ onto its Euclidean counterpart on $\mathbb{R}^N$ while preserving the $L^p$ norm, together with a one-dimensional tangent line inequality, we establish such sharp Poincar\'e--Sobolev inequalities (namely, Sobolev inequalities with the sharp Euclidean constant after the full spectral gap is subtracted) for $N=3$ and $2\leq p<3$, for $N=4$ and $\frac{1+\sqrt{13}}{2}\leq p<4$, and for $N\geq5$ and $\min\{3,\frac{1+\sqrt{4N-3}}{2}\}\leq p<N$. In particular,  in all cases we also obtain an explicit remainder term. Moreover, since the change of function keeps the $L^p$ norm, we prove the sharp logarithmic Sobolev inequalities with the full spectral gap and the sharp Euclidean constant for $N=3$ and $2\leq p\leq3$, and for $N\geq4$ on an explicit interval of $p$. For $p\geq3$, we also use the change of variable $t=\sinh\rho$ and the Hardy inequality, which give a second remainder term. At the endpoint $p=N$, this remainder gives a Hardy--Trudinger--Moser inequality on $\mathbb{H}^N$ under a weaker constraint than the one of Liang, Lu, Wang and Yang, in every dimension $N\geq3$. Finally, we show that for $1<p<\frac{1+\sqrt{4N-3}}{2}$, which includes all $1<p<2$, the full spectral gap can still be subtracted with some positive constant in front of the Sobolev term, and that this constant is smaller than $S(N,p)^{p}$ when $p<\min\{2,\sqrt N\}$. Together with our other results, this shows that the full spectral gap can be subtracted for every $1<p<N$.
\end{abstract}
\maketitle

\section{Introduction}\label{secintro}

The main goal of this paper is to study sharp Poincar\'e--Sobolev inequalities and sharp logarithmic Sobolev inequalities on the hyperbolic space. We want to know when the sharp Sobolev constant of the Euclidean space $\mathbb{R}^N$ still works on the hyperbolic space after we subtract the full spectral gap of the $p$-Laplacian from the Dirichlet integral. We also want to know if the inequality can be improved by a remainder term.

In this paper we work with the Poincar\'e ball model of the hyperbolic space $\mathbb{H}^N$, $N\geq3$. In this model, $\mathbb{H}^N$ is the unit ball
\[
B^N := \{ x \in \mathbb{R}^N : |x| < 1 \}
\]
with the Riemannian metric
\[
g(x) = \left( \frac{2}{1 - |x|^2} \right)^2 dx \otimes dx .
\]
It is a complete, non-compact Riemannian manifold with constant sectional curvature $-1$, see \cite{RAT}. The volume element is
\[
dV_{\mathbb{H}} = \left( \frac{2}{1 - |x|^{2}} \right)^{N} dx ,
\]
and the gradient is
\[
\nabla_{\mathbb{H}} = \left( \frac{1 - |x|^{2}}{2} \right)^{2} \nabla ,
\]
where $\nabla$ is the Euclidean gradient. The hyperbolic distance from a point $x$ to the origin is
\[
\rho(x) = \int_{0}^{|x|} \frac{2}{1-s^{2}}\,ds = \ln\frac{1+|x|}{1-|x|},
\qquad\text{so that}\qquad |x| = \tanh\frac{\rho(x)}{2}.
\]
In geodesic polar coordinates about the origin we have
\[
dV_{\mathbb{H}}=\sinh^{N-1}\rho\,d\rho\,d\sigma ,
\]
where $d\sigma$ is the surface measure of the unit sphere $\mathbb{S}^{N-1}$. We write $|\mathbb{S}^{N-1}|$ for the area of $\mathbb{S}^{N-1}$. If $u=u(\rho)$ is a radial function, then $|\nabla_{\mathbb{H}}u|=|u'(\rho)|$. We denote by $B_R$ the geodesic ball of radius $R$ with center at the origin. Finally, $D^{1,p}(\mathbb{H}^N)$ is the completion of $C_c^\infty(\mathbb{H}^N)$ with respect to the norm $\|\nabla_{\mathbb{H}}u\|_{L^p(\mathbb{H}^N)}$.

\medskip
\subsection{Sobolev and Poincar\'e inequalities on the hyperbolic space.}

On the hyperbolic space $\mathbb{H}^{N}$, the bottom of the spectrum of the operator
\[
-\operatorname{div}\left(|\nabla_{\mathbb{H}} u|^{p-2}\nabla_{\mathbb{H}} u\right), \qquad p>1,
\]
is positive. This is different from the Euclidean space, where the spectrum of the $p$-Laplacian is $[0,\infty)$. More precisely, for every $p>1$, we have the Poincar\'e inequality
\begin{equation}\label{eq:gap}
\int_{\mathbb{H}^N} |\nabla_{\mathbb{H}} u|^p \, dV_{\mathbb{H}}
\geq \left(\frac{N-1}{p}\right)^p
\int_{\mathbb{H}^N} |u|^p \, dV_{\mathbb{H}},
\qquad u \in C_c^\infty(\mathbb{H}^N).
\end{equation}
The constant $\left(\frac{N-1}{p}\right)^p$ is sharp and it is never attained. See \cite{BG16, BGG17, FLL22, FLLM23, KS16, MS08, NN19, H21}, for instance.

For $1\leq p<N$, let $S(N,p)$ be the sharp constant in the $L^p$ Sobolev inequality on $\mathbb{R}^N$,
\begin{equation}\label{eq:sobolev}
\int_{\mathbb{R}^N} |\nabla w|^p \, dx
\geq S(N,p)^p
\left( \int_{\mathbb{R}^N} |w|^{\frac{Np}{N-p}} \, dx \right)^{\frac{N-p}{N}},
\qquad w\in C_c^\infty(\mathbb{R}^N).
\end{equation}
We refer to Aubin \cite{Aub76} and Talenti \cite{T76} for the value of $S(N,p)$ and for the extremal functions. The same inequality holds on $\mathbb{H}^N$ with the same constant, see \cite{ADH98}.

When $p=2$ and $N\geq3$, one has the Poincar\'e--Sobolev inequality
\begin{equation}\label{PoinSob}
\int_{\mathbb{H}^N} |\nabla_{\mathbb{H}} u|^2 \, dV_{\mathbb{H}}
- \frac{(N-1)^2}{4} \int_{\mathbb{H}^N} |u|^2 \, dV_{\mathbb{H}}
\geq C(N)^2
\left( \int_{\mathbb{H}^N} |u|^{\frac{2N}{N-2}} \, dV_{\mathbb{H}} \right)^{\frac{N-2}{N}},
\qquad u \in C_c^\infty(\mathbb{H}^N),
\end{equation}
for some constant $C(N)>0$, see \cite[Appendix B]{MS08}. In dimension three, Benguria, Frank and Loss \cite{BFL08} proved that the best constant in \eqref{PoinSob} is the sharp Euclidean constant, that is,
\[
C(3)=S(3,2).
\]
In their work, \eqref{PoinSob} with $N=3$ is the same as the sharp Hardy--Sobolev--Maz'ya inequality on the half space. Mancini and Sandeep \cite{MS08} got the same result from the fact that the critical equation on $\mathbb{H}^3$ has no positive solution. Lu and Yang \cite{LY19} gave a third proof by developing a method of using the Helgason--Fourier analysis on $\mathbb{H}^N$ to establish sharp geometric inequalities. Their method also works for higher order operators. Through such an approach,  they proved Poincar\'e--Sobolev inequalities for the Paneitz and GJMS operators on $\mathbb{H}^{N}$, and the higher order Hardy--Sobolev--Maz'ya inequalities on half spaces. On the other hand, Hebey \cite{Heb99} proved that
\[
C(N)<S(N,2)\qquad\text{for }N\geq4 .
\]
Tertikas and Tintarev \cite{TT07} proved that \eqref{PoinSob} has extremal functions when $N\geq 4$. The inequality with the constant $S(N,2)$ holds if we subtract the smaller number $\frac{N(N-2)}{4}$ in place of $\frac{(N-1)^{2}}{4}$, see \cite{Heb99}. The higher order picture is the same.
  Lu and Yang \cite{LY221} proved that for the $k$-th order inequality the best constant is the Sobolev constant when $N=2k+1$, and is strictly smaller when $N\geq2k+2$. For the complex, quaternionic and Cayley hyperbolic spaces, we refer to \cite{FLY24, LY22}. Lu and Tao \cite{LuTao} recently  proved the existence of extremal functions for the higher order Poincar\'e--Sobolev and Hardy--Sobolev--Maz'ya inequalities when $k\geq2$ and $N\geq2k+2$, and from this the existence of positive solutions of the higher order Brezis--Nirenberg equation for the GJMS operator $P_k$ on $\mathbb{H}^N$:
\begin{equation}\label{sublambdaBN}
P_{k}\left(f\right)-\alpha f=|f|^{p-2}f
\end{equation}
for $\alpha=\prod\limits_{i=1}^{k}\frac{\left(2i-1\right)^2}{4}$. Earlier works on $\mathbb{H}^n$ include: 
For $k=1$, Mancini and Sandeep \cite{MS08} showed that an entire solution exists on $\mathbb{H}^n$ either when  $n\geq4$, $p=\frac{2n}{n-2}$, $0<\alpha\leq\frac{1}{4}$ or when $n\geq3$, $1<p<\frac{2n}{n-2}$ and $\alpha\leq\frac{1}{4}$.  
 For $k\geq 2$, Li et al. \cite{LLY-AIM} established the existence of a nontrivial solution  when $n\geq 2k+2$, $p=\frac{2n}{n-2k}$, $\alpha<\prod\limits_{i=1}^{k}\frac{\left(2i-1\right)^2}{4}$.

\medskip
\subsection{$AB$ program.}
The $AB$ program of Aubin \cite{Aub76}, which was developed by Hebey \cite{Heb99} and by Druet and Hebey \cite{DH02}, is about inequalities of the form
\begin{equation}\tag{$AB_{p}$}
A\left(\int_{\mathbb{H}^N}|u|^{\frac{Np}{N-p}}dV_{\mathbb{H}}\right)^{\frac{N-p}{N}}\leq \int_{\mathbb{H}^N}|\nabla_{\mathbb{H}}u|^{p}dV_{\mathbb{H}}+B\int_{\mathbb{H}^N}|u|^{p}dV_{\mathbb{H}},\qquad u\in C_{c}^{\infty}(\mathbb{H}^N).
\end{equation}
The best first constant is the largest $A$ such that $(AB_{p})$ holds with some $B$. The best second constant is the smallest $B$ such that $(AB_{p})$ holds with some $A$. The $A$-part of the program asks if the best first constant is $S(N,p)^{p}$, and if $(AB_{p})$ holds with $A=S(N,p)^{p}$. On Cartan--Hadamard manifolds which satisfy the Cartan--Hadamard conjecture, and in particular on $\mathbb{H}^N$, this was proved by Aubin, Druet and Hebey \cite{ADH98} with $B=0$. The $B$-part of the program asks for the best second constant. By \eqref{eq:gap} and the sharpness of $\left(\frac{N-1}{p}\right)^{p}$, we cannot take $B<-\left(\frac{N-1}{p}\right)^{p}$. So the best second constant is $-\left(\frac{N-1}{p}\right)^{p}$ if $(AB_{p})$ holds with $B=-\left(\frac{N-1}{p}\right)^{p}$ and some $A>0$. This is true for $p=2$ by \eqref{PoinSob}. We prove in Corollary \ref{cor:B} below that it is true for every $1<p<N$.

The question that is left is if the two best constants can be taken at the same time. In the language of \cite{DH02}, this is the question of the second best constant when the first constant is sharp. In other words, we ask if
\begin{equation}\label{eq:PS}
\int_{\mathbb{H}^N}|\nabla_{\mathbb{H}}u|^{p}dV_{\mathbb{H}}-\left(\frac{N-1}{p}\right)^{p}\int_{\mathbb{H}^N}|u|^{p}dV_{\mathbb{H}}
\geq S(N,p)^{p}\left(\int_{\mathbb{H}^N}|u|^{\frac{Np}{N-p}}dV_{\mathbb{H}}\right)^{\frac{N-p}{N}}
\end{equation}
holds for all $u\in C_{c}^{\infty}(\mathbb{H}^N)$. When it holds, both constants in \eqref{eq:PS} are sharp, and we call \eqref{eq:PS} a doubly sharp inequality. By the results above, for $p=2$ the inequality \eqref{eq:PS} holds when $N=3$ and it fails when $N\geq4$.

For $p\neq 2$ much less is known. Nguyen \cite{H18} proved \eqref{eq:PS} for
\begin{equation}\label{eq:Nguyen}
N\geq4\qquad\text{and}\qquad\frac{2N}{N-1}\leq p<N ,
\end{equation}
together with sharp Poincar\'e--Gagliardo--Nirenberg inequalities. His method is to compare the Dirichlet integrals of the hyperbolic and the Euclidean symmetric decreasing rearrangements of $u$. At the endpoint $p=N$, Liang, Lu, Wang and Yang \cite{LLWY} proved the Hardy--Trudinger--Moser inequality on $\mathbb{H}^N$, in which an exponential integral is bounded when the left side of \eqref{eq:PS} with $p=N$ is at most $1$. In \cite{DGLL}, we proved that for every $N\geq3$ and $2\leq p<N$ there is a constant $\lambda(N,p)>0$ such that
\[
\int_{\mathbb{H}^N} |\nabla_{\mathbb{H}} u|^p \, dV_{\mathbb{H}}
- \lambda(N,p) \int_{\mathbb{H}^N} |u|^p \, dV_{\mathbb{H}}
\geq S(N,p)^p
\left( \int_{\mathbb{H}^N} |u|^{\frac{Np}{N-p}} \, dV_{\mathbb{H}} \right)^{\frac{N-p}{N}} ,
\qquad u \in C_c^\infty(\mathbb{H}^N).
\]
But $\lambda(N,p)<\left(\frac{N-1}{p}\right)^{p}$ outside the range \eqref{eq:Nguyen}. So the case $N=3$, which contains the inequality of Benguria, Frank and Loss, and the range $2<p<\frac{2N}{N-1}$ for $N\geq4$ were still open.

\medskip
\subsection{Logarithmic Sobolev inequalities.}
The same question can be asked for the logarithmic Sobolev inequality. The logarithmic Sobolev inequality goes back to Stam \cite{Stam59} and Federbush \cite{Fede69} in information theory. Gross \cite{Gro75} proved it for the Gaussian measure:
\begin{equation}\label{eq:gross}
\int_{\mathbb{R}^N} |u|^2 \ln |u|^2 \, d\mu
\leq 2 \int_{\mathbb{R}^N} |\nabla u|^2 \, d\mu ,
\qquad\text{if}\qquad \int_{\mathbb{R}^N} |u|^2 \, d\mu = 1,
\end{equation}
where
\[
d\mu(x) = (2\pi)^{-\frac{N}{2}} e^{-\frac{|x|^2}{2}} \, dx
\]
is the standard Gaussian measure on $\mathbb{R}^N$. The constant $2$ in \eqref{eq:gross} is sharp, and equality holds for $u(x)=A e^{b\cdot x}$. By a change of variables, one can get from \eqref{eq:gross} the logarithmic Sobolev inequality for the Lebesgue measure, see \cite{Wei78},
\begin{equation}\label{eq:logE2}
\int_{\mathbb{R}^N} |u|^2 \ln |u| \, dx
\leq \frac{N}{4} \ln \left( \frac{2}{\pi N e}
\int_{\mathbb{R}^N} |\nabla u|^2 \, dx \right),
\qquad\text{if}\qquad \int_{\mathbb{R}^N} |u|^2 \, dx = 1.
\end{equation}
Logarithmic Sobolev inequalities are now basic tools in analysis and probability, see \cite{BGL14} for an overview and \cite{WB-99, YF, IG, LL18} for some of the many later works. The sharp $L^p$ version of \eqref{eq:logE2} holds for every $p>1$ in every dimension, see \cite{PD03, YF, IG}. It says that
\begin{equation}\label{eq:logE}
\int_{\mathbb{R}^N}|v|^p\ln|v|\,dx \leq \frac{N}{p^2}\ln \left[\mathcal{L}_{N,p} \int_{\mathbb{R}^N}|\nabla v|^pdx\right],
\qquad\text{if}\qquad \int_{\mathbb{R}^N}|v|^pdx=1,
\end{equation}
where the sharp constant is
\begin{equation}\label{eq:LNp}
\mathcal{L}_{N,p} = \frac{p}{N}\pi^{-\frac{p}{2}} \left( \frac{p-1}{e} \right)^{p-1} \left( \frac{\Gamma\left( \frac{N}{2} + 1 \right) }{\Gamma\left( \frac{N(p-1)}{p} + 1 \right)} \right)^{\frac{p}{N}}.
\end{equation}
For $p=2$ we have $\mathcal{L}_{N,2}=\frac{2}{\pi Ne}$, which is the constant in \eqref{eq:logE2}.

On Riemannian manifolds, the logarithmic Sobolev constant depends on the geometry of the manifold, for example on the dimension, the diameter and lower bounds of the Ricci curvature. Cavalletti and Mondino \cite{CM} proved sharp Sobolev and logarithmic Sobolev inequalities on spaces with the condition $CD(K,N)$, $K>0$, by optimal transport. For noncompact and nonconvex manifolds, we refer to the works of F.-Y. Wang \cite{FW-97, FW-99, FW-09, FW-209}.

On the hyperbolic space the curvature is negative. In the spirit of the $AB$ program, we ask if the following inequality holds:
\begin{equation}\label{eq:logPS}
\int_{\mathbb{H}^N} |u|^p \ln |u| \, dV_{\mathbb{H}}
\leq \frac{N}{p^2}
\ln \left[
\mathcal{L}_{N,p}
\left(
\int_{\mathbb{H}^N} |\nabla_{\mathbb{H}} u|^p \, dV_{\mathbb{H}} - \left(\frac{N-1}{p}\right)^{p}
\right)
\right],
\qquad
\int_{\mathbb{H}^N} |u|^p \, dV_{\mathbb{H}} = 1 .
\end{equation}
As for \eqref{eq:PS}, the constant $\left(\frac{N-1}{p}\right)^{p}$ cannot be replaced by a larger one. In \cite{H18}, Nguyen proved \eqref{eq:logPS} in the range \eqref{eq:Nguyen}, as a limit of his Poincar\'e--Gagliardo--Nirenberg inequalities. His comparison of the hyperbolic and the Euclidean rearrangements \cite[Theorem 2.2]{H18} holds for every $p\geq\frac{2N}{N-1}$, with no upper bound on $p$. So, together with \eqref{eq:logE} for all $p>1$ from \cite{IG}, it gives \eqref{eq:logPS} for all $N\geq3$ and $p\geq\frac{2N}{N-1}$. In \cite{DGLL}, we proved \eqref{eq:logPS} for all $N\geq3$ and $2\leq p<N$ with $\lambda(N,p)$ in place of $\left(\frac{N-1}{p}\right)^{p}$, with the explicit value $\lambda(N,2)=\frac{N^{2}(N-1)}{4(N+2)}$. We also used it to get logarithmic Sobolev and Beckner inequalities for Gaussian type measures on $\mathbb{H}^N$. Again the case $N=3$ and $2\leq p<3$ was open, and the case $p=2$ was open in every dimension.

\medskip
\subsection{Main results.}\label{sec:results}
We state our results on \eqref{eq:PS} in two forms. The first form is the inequality \eqref{eq:PS} without remainder. Here we give the largest range of $p$ that our methods can reach. The second form has an explicit remainder term. There are two remainders. The first one is a weight $R_{N}$, which we have for $p$ up to $3$, and for $N\geq6$ also for some $p>3$. The second one is an exponentially small weight, which we have for all $p\geq3$. The range for the first remainder is smaller than the range of the inequality without remainder, because we need $R_{N}$ to be monotone when we pass from radial functions to all functions.

Our first main result is the following doubly sharp Poincar\'e--Sobolev inequality.

\begin{theorem}\label{thm:PS}
Let $N\geq3$ and assume that one of the following conditions holds:
\begin{itemize}
\item[(i)] $N=3$ and $2\leq p<3$,
\item[(ii)] $N=4$ and $\frac{1+\sqrt{13}}{2}\leq p<4$,
\item[(iii)] $N\geq5$ and $\min\left\{3,\frac{1+\sqrt{4N-3}}{2}\right\}\leq p<N$.
\end{itemize}
Then for all $u\in D^{1,p}(\mathbb{H}^N)$, we have
\[
\int_{\mathbb{H}^N}|\nabla_{\mathbb{H}}u|^{p}dV_{\mathbb{H}}-\left(\frac{N-1}{p}\right)^{p}\int_{\mathbb{H}^N}|u|^{p}dV_{\mathbb{H}}
\geq S(N,p)^{p}\left(\int_{\mathbb{H}^N}|u|^{\frac{Np}{N-p}}dV_{\mathbb{H}}\right)^{\frac{N-p}{N}}.
\]
Both constants $\left(\frac{N-1}{p}\right)^{p}$ and $S(N,p)^{p}$ are sharp, and $S(N,p)^{p}$ is not attained.
\end{theorem}

\begin{remark}
{\rm
When $N=3$ and $p=2$, Theorem \ref{thm:PS} is the inequality of Benguria, Frank and Loss. The number $\frac{1+\sqrt{4N-3}}{2}$ is the positive root of $p^{2}-p=N-1$, that is, of $p-1=\frac{N-1}{p}$. It is larger than $2$ for $N\geq4$, which agrees with the theorem of Hebey, and it is equal to $3$ for $N=7$. For $N=4$ it is $\frac{1+\sqrt{13}}{2}<\frac83$. So Theorem \ref{thm:PS} is new for $N=3$ and $2<p<3$, and for $N=4$ and $\frac{1+\sqrt{13}}{2}\leq p<\frac83$. For $N\geq5$, the range in (iii) is contained in the range \eqref{eq:Nguyen} of Nguyen, since $3\geq\frac{2N}{N-1}$ and $\frac{1+\sqrt{4N-3}}{2}>\frac{2N}{N-1}$ for $N\geq5$. In this case the new results are the remainder terms below. We will also show in Proposition \ref{prop:fail} that \eqref{eq:PS} fails for $p=1$ and for $1<p<\min\{2,\sqrt N\}$. So for $N=3$ the lower bound $p=2$ in (i) is close to the best one.
}
\end{remark}

Our first remainder term uses the function
\begin{equation}\label{eq:beta}
R_{N}(r):=\left(p-1-\frac{N-1}{p}\right)\left[\frac{1}{r^{p}}-\frac{\coth^{p-2}r}{\sinh^{2}r}\right]
+\frac{N-1}{p}\left(\coth^{p-2}r-1\right),\qquad r>0 .
\end{equation}
For $N=3$ it is
\[
R_{3}(r)=\frac{(p-2)(p+1)}{p}\left[\frac{1}{r^{p}}-\frac{\coth^{p-2}r}{\sinh^{2}r}\right]+\frac2p\left(\coth^{p-2}r-1\right).
\]

\begin{theorem}\label{thm:main}
Let $N\geq3$ and assume that one of the following conditions holds:
\begin{itemize}
\item[(i)] $N=3$ and $2\leq p<3$,
\item[(ii)] $N=4$ and $\frac{1+\sqrt{13}}{2}\leq p\leq3$,
\item[(iii)] $N\geq5$ and $\frac{1+\sqrt{4N-3}}{2}\leq p\leq\min\left\{5,\frac{N+1}{2}\right\}$
(this range is not empty for $5\leq N\leq21$).
\end{itemize}
Then for all $u\in D^{1,p}(\mathbb{H}^N)$, we have
\begin{equation}\label{eq:main}
\begin{aligned}
\int_{\mathbb{H}^N}|\nabla_{\mathbb{H}}u|^{p}dV_{\mathbb{H}}-\left(\frac{N-1}{p}\right)^{p}\int_{\mathbb{H}^N}|u|^{p}dV_{\mathbb{H}}
\geq\ &S(N,p)^{p}\left(\int_{\mathbb{H}^N}|u|^{\frac{Np}{N-p}}dV_{\mathbb{H}}\right)^{\frac{N-p}{N}}\\
&+\left(\frac{N-1}{p}\right)^{p-1}\int_{\mathbb{H}^N}R_{N}(\rho)\,|u|^{p}dV_{\mathbb{H}} .
\end{aligned}
\end{equation}
Moreover, the function $R_{N}$ is nonincreasing on $(0,\infty)$ and $R_{N}(r)\to0$ as $r\to\infty$. In particular, $R_{N}\geq0$. When $N=3$ and $p=2$, we have $R_{3}\equiv0$. In all other cases of (i)--(iii), $R_{N}$ is strictly decreasing and strictly positive, and
\begin{equation}\label{eq:Basym}
\lim_{r\to0^{+}}r^{p-2}R_{N}(r)=\frac{N-(p-2)^{2}}{3}\geq0,\qquad
\lim_{r\to\infty}r^{p}R_{N}(r)=p-1-\frac{N-1}{p}\geq0 .
\end{equation}
\end{theorem}
\begin{remark}
{\rm
The two limits in \eqref{eq:Basym} explain the range of Theorem \ref{thm:main}. The limit at infinity is nonnegative if and only if $p-1\geq\frac{N-1}{p}$, which is the lower bound in (ii) and (iii). The limit at the origin is positive if and only if $(p-2)^{2}<N$, that is, $p<2+\sqrt N$. So $R_{N}$ can be nonnegative only when $\frac{1+\sqrt{4N-3}}{2}\leq p\leq2+\sqrt N$. The upper bounds $p\leq3$ in (ii) and $p\leq\min\{5,\frac{N+1}{2}\}$ in (iii) come from our proof of the monotonicity.
}
\end{remark}
When $3\leq p<N$, we have the second remainder term:

\begin{theorem}\label{thm:pge3}
Let $N\geq4$ and $3\leq p<N$. Then for all $u\in D^{1,p}(\mathbb{H}^N)$, we have
\begin{equation}\label{eq:mainW}
\begin{aligned}
\int_{\mathbb{H}^N}|\nabla_{\mathbb{H}}u|^{p}dV_{\mathbb{H}}-\left(\frac{N-1}{p}\right)^{p}\int_{\mathbb{H}^N}|u|^{p}dV_{\mathbb{H}}
\geq\ &S(N,p)^{p}\left(\int_{\mathbb{H}^N}|u|^{\frac{Np}{N-p}}dV_{\mathbb{H}}\right)^{\frac{N-p}{N}}\\
&+\left(\frac{N-1}{p}\right)^{p}\int_{\mathbb{H}^N}\frac{2}{e^{2\rho}-1}\,|u|^{p}dV_{\mathbb{H}} .
\end{aligned}
\end{equation}
\end{theorem}

The proof of Theorem \ref{thm:pge3} is different from the proof of Theorem \ref{thm:main}. We change the radial variable to $t=\sinh\rho$, as in \cite{LLWY} for $p=N=3$, and we use the one-dimensional Hardy inequality. This gives the inequality in every dimension with a very short proof. The weight $\frac{2}{e^{2\rho}-1}$ is decreasing. It behaves like $\rho^{-1}$ near the origin and it is exponentially small at infinity. Note that $\frac{2}{e^{2\rho}-1}=\coth\rho-1$. So for $p=3$, by the definition \eqref{eq:beta} of $R_N$,
\[
\left(\frac{N-1}{3}\right)^{2}R_N(\rho)-\left(\frac{N-1}{3}\right)^{3}\frac{2}{e^{2\rho}-1}=\left(\frac{N-1}{3}\right)^{2}\left(2-\frac{N-1}{3}\right)\left[\frac{1}{\rho^{3}}-\frac{\coth\rho}{\sinh^{2}\rho}\right]\geq0
\]
for $4\leq N\leq7$, which are the dimensions where both theorems hold for $p=3$. Here we used $N\leq7$ and the inequality \eqref{eq:Laz} of Lazarevi\'c below. So for $p=3$ the remainder in Theorem \ref{thm:main} is larger than the remainder in Theorem \ref{thm:pge3}. For $p>3$, in the cases where both theorems hold, the weight $R_{N}$ is also larger near the origin, where it behaves like $\rho^{2-p}$, and at infinity, where it behaves like $\left(p-1-\frac{N-1}{p}\right)\rho^{-p}$. The advantage of Theorem \ref{thm:pge3} is that it holds for all $3\leq p<N$. Theorem \ref{thm:PS} follows from Theorems \ref{thm:main} and \ref{thm:pge3} if we drop the remainder terms, since the union of the ranges of Theorems \ref{thm:main} and \ref{thm:pge3} is the range of Theorem \ref{thm:PS}.

Below the range of Theorem \ref{thm:PS}, we can still subtract the full spectral gap, but the constant in front of the Sobolev term is smaller. This includes all $1<p<2$.

\begin{theorem}\label{thm:gap}
Let $N\geq3$ and $1<p<\frac{1+\sqrt{4N-3}}{2}$. Then there is a constant $C=C(N,p)>0$ such that for all $u\in D^{1,p}(\mathbb{H}^N)$, we have
\begin{equation}\label{eq:gapC}
\int_{\mathbb{H}^N}|\nabla_{\mathbb{H}}u|^{p}dV_{\mathbb{H}}-\left(\frac{N-1}{p}\right)^{p}\int_{\mathbb{H}^N}|u|^{p}dV_{\mathbb{H}}
\geq C\left(\int_{\mathbb{H}^N}|u|^{\frac{Np}{N-p}}dV_{\mathbb{H}}\right)^{\frac{N-p}{N}}.
\end{equation}
The best constant $C$ in \eqref{eq:gapC} is at most $S(N,p)^{p}$, and it is strictly smaller than $S(N,p)^{p}$ when $p<\min\{2,\sqrt N\}$.
\end{theorem}

Note that $\frac{1+\sqrt{4N-3}}{2}<N$ for $N\geq3$, so $p<N$ in Theorem \ref{thm:gap}. For $N=3$ the range of Theorem \ref{thm:gap} is $1<p<2$, and for $N\geq4$ it contains $1<p\leq2$. For $p=2$ and $N\geq4$, Theorem \ref{thm:gap} gives another proof of \eqref{PoinSob}. The upper bound $S(N,p)^{p}$ for the best constant follows from the test functions in the proof of the sharpness of $S(N,p)^{p}$ in Section \ref{sec:proof}, and the strict inequality follows from Proposition \ref{prop:fail}. In the proof of Theorem \ref{thm:gap} we use the ground state transform $u=g\sinh^{-\frac{N-1}{p}}\rho$ and a convexity estimate for the function $t\mapsto|t|^{p}$. For $1<p<2$ this estimate is of the type used by Banerjee, Kijaczko and Sahu \cite{BKS26} for Hardy--Sobolev--Maz'ya inequalities on $\mathbb{R}^N$ with $1<p<2$. Since we do not look for the sharp constant in \eqref{eq:gapC}, we do not need the change of function \eqref{eq:mapintro} and the function $R_{N}$.

If we put Theorem \ref{thm:PS} and Theorem \ref{thm:gap} together, we get the $B$-part of the $AB$ program on $\mathbb{H}^N$ for all $p$.

\begin{corollary}\label{cor:B}
Let $N\geq3$ and $1<p<N$. Then there is a constant $C=C(N,p)>0$ such that \eqref{eq:gapC} holds for all $u\in D^{1,p}(\mathbb{H}^N)$. So the best second constant in $(AB_{p})$ is $-\left(\frac{N-1}{p}\right)^{p}$ for every $1<p<N$.
\end{corollary}

Indeed, for $p<\frac{1+\sqrt{4N-3}}{2}$ this is Theorem \ref{thm:gap}. For $p\geq\frac{1+\sqrt{4N-3}}{2}$ it follows from Theorem \ref{thm:PS}, since the lower bound for $p$ in each of the cases (i)--(iii) of Theorem \ref{thm:PS} is at most $\frac{1+\sqrt{4N-3}}{2}$. For $N=3$ this lower bound is $2=\frac{1+\sqrt{9}}{2}$.

Our second main result is about the logarithmic Sobolev inequality \eqref{eq:logPS}. The change of function in the proof of Theorem \ref{thm:main} keeps the $L^{p}$ norm. This is exactly what the logarithmic Sobolev inequality needs, and here we do not need the condition $p<N$. Again we state the result in two forms.

\begin{theorem}\label{thm:logPS}
Let $N\geq3$ and $p\geq\frac{1+\sqrt{4N-3}}{2}$ (note that $\frac{1+\sqrt{4N-3}}{2}\geq2$). Assume that one of the following conditions holds:
\begin{itemize}
\item[(a)] $p\leq3$,
\item[(b)] $p\leq\min\left\{5,\frac{N+1}{2}\right\}$,
\item[(c)] $5\leq p\leq\frac{3+\sqrt{4N-3}}{2}$.
\end{itemize}
Then for all $u\in W^{1,p}(\mathbb{H}^N)$ with $\int_{\mathbb{H}^N}|u|^{p}dV_{\mathbb{H}}=1$, we have
\[
\int_{\mathbb{H}^N} |u|^p \ln |u| \, dV_{\mathbb{H}}
\leq \frac{N}{p^2}
\ln \left[
\mathcal{L}_{N,p}
\left(
\int_{\mathbb{H}^N} |\nabla_{\mathbb{H}} u|^p \, dV_{\mathbb{H}} - \left(\frac{N-1}{p}\right)^{p}
\right)
\right].
\]
Both constants $\left(\frac{N-1}{p}\right)^{p}$ and $\mathcal{L}_{N,p}$ are sharp, and $\mathcal{L}_{N,p}$ is not attained.
\end{theorem}

For each $N$, the range of $p$ in Theorem \ref{thm:logPS} is as follows: $2\leq p\leq3$ for $N=3$, $\frac{1+\sqrt{4N-3}}{2}\leq p\leq3$ for $N=4,5$, $\frac{1+\sqrt{4N-3}}{2}\leq p\leq\frac{N+1}{2}$ for $N=6,7,8$, $\frac{1+\sqrt{4N-3}}{2}\leq p\leq5$ for $9\leq N\leq13$, and $\frac{1+\sqrt{4N-3}}{2}\leq p\leq\frac{3+\sqrt{4N-3}}{2}$ for $N\geq14$. For $N=3$, Theorem \ref{thm:logPS} gives the logarithmic Sobolev inequality with the full spectral gap for all $2\leq p\leq3$. This includes the case $p=2$, which was open, and the endpoint $p=N=3$. For $N=4$, it gives the inequality for $\frac{1+\sqrt{13}}{2}\leq p\leq3$, which includes the range $p<\frac83$ that is not covered by \eqref{eq:Nguyen}.

If we combine Theorem \ref{thm:logPS} with the result of Nguyen \cite[Theorem 2.2]{H18} and the argument of \cite[Section 3]{DGLL}, we get the following. Note that $\frac{2N}{N-1}\leq3$ for $N\geq3$ and $\frac{2N}{N-1}<\frac{1+\sqrt{4N-3}}{2}$ for $N\geq5$.

\begin{corollary}\label{cor:log}
The inequality \eqref{eq:logPS} holds for all $p\geq2$ when $N=3$, for all $p\geq\frac{1+\sqrt{13}}{2}$ when $N=4$, and for all $p\geq\frac{2N}{N-1}$ when $N\geq5$.
\end{corollary}

The logarithmic Sobolev inequality with remainder is the following.

\begin{theorem}\label{thm:log}
Let $N\geq3$ and assume that either $N=3$ and $2\leq p\leq3$, or $(N,p)$ satisfies (ii) or (iii) of Theorem \ref{thm:main}. Then for all $u\in W^{1,p}(\mathbb{H}^N)$ with $\int_{\mathbb{H}^N}|u|^{p}dV_{\mathbb{H}}=1$, we have
\begin{equation}\label{eq:log}
\begin{aligned}
\int_{\mathbb{H}^N}|u|^{p}\ln|u|\,dV_{\mathbb{H}}\leq \frac{N}{p^{2}}\ln\Bigg[\mathcal{L}_{N,p}\Bigg(&
\int_{\mathbb{H}^N}|\nabla_{\mathbb{H}}u|^{p}dV_{\mathbb{H}}-\left(\frac{N-1}{p}\right)^{p}\\
&-\left(\frac{N-1}{p}\right)^{p-1}\int_{\mathbb{H}^N}R_{N}(\rho)\,|u|^{p}dV_{\mathbb{H}}\Bigg)\Bigg].
\end{aligned}
\end{equation}
Moreover, if $u$ is radial and nonincreasing in $\rho$, then the term
\[
\frac{N-1}{p}\int_{\mathbb{H}^N}|u|^{p}\ln\frac{\sinh\rho}{\rho}\,dV_{\mathbb{H}}
\]
can be added to the left side of \eqref{eq:log}.
\end{theorem}

The remainder inside the logarithm in \eqref{eq:log} is the same as in Theorem \ref{thm:main}. The range is also the same as in Theorem \ref{thm:main}, because we need the monotonicity of $R_{N}$ in the Hardy--Littlewood inequality. The only difference is that the endpoint $p=N=3$ is allowed, since the condition $p<N$ is not needed here. The second remainder for radial functions comes from the change of function, and it is positive except at the origin. The change of variable $t=\sinh\rho$ in the proof of Theorem \ref{thm:pge3} does not keep the $L^{p}$ norm. So it does not give a logarithmic Sobolev inequality.

Finally, we study the endpoint $p=N$. For $N\geq2$, let
\[
\Phi_{N}(s)=e^{s}-\sum_{k=0}^{N-1}\frac{s^{k}}{k!},\qquad
\beta_{N}=N\left(\frac{N\pi^{N/2}}{\Gamma(N/2+1)}\right)^{\frac{1}{N-1}} .
\]
Here $\beta_{N}$ is the Moser constant. Trudinger--Moser and Adams inequalities on hyperbolic spaces have been studied by  many authors.  Mancini, Sandeep and Tintarev \cite{MST13, KS16},  Lu and Tang \cite{LT13, LT16} found the best constants in the Trudinger--Moser inequalities on hyperbolic spaces of any dimension, and also in the form with the exact growth condition in \cite{LT16}.  Lu and Yang \cite{LY17} proved the sharp Hardy--Adams inequality for the bi-Laplacian on $\mathbb{H}^4$, and Li, Lu and Yang extended it to hyperbolic spaces of any even dimension \cite{LLY18}, and then to Adams and Hardy--Adams inequalities of any fractional order on hyperbolic spaces of all dimensions \cite{LLY20}. Lu and Yang \cite{LY22} and Flynn, Lu and Yang \cite{FLY24} proved sharp Adams and Hardy--Adams inequalities on the complex and quaternionic hyperbolic spaces and on the Cayley hyperbolic plane.  do \'O, Lu and Ponciano \cite{DLP25, DLP26} proved sharp Adams--Trudinger--Moser inequalities for symmetric functions without boundary conditions on hyperbolic spaces and spherically symmetric Riemannian manifolds. The main result of Liang et al \cite{LLWY} is that for $N\geq3$ there is a constant $C=C(N)$ such that
\begin{equation}\label{eq:HTM}
\int_{\mathbb{H}^N}\Phi_{N}\left(\beta_{N}|u|^{\frac{N}{N-1}}\right)dV_{\mathbb{H}}\leq C
\end{equation}
for all $u\in W^{1,N}_{0}(\mathbb{H}^N)$ with
\[
\int_{\mathbb{H}^N}|\nabla_{\mathbb{H}}u|^{N}dV_{\mathbb{H}}-\left(\frac{N-1}{N}\right)^{N}\int_{\mathbb{H}^N}|u|^{N}dV_{\mathbb{H}}\leq1 .
\]
The proof in \cite{LLWY} uses a lower bound of the left side of the constraint by a Dirichlet integral on $\mathbb{R}^N$ plus a weighted $L^{N}$ norm, see \cite[Lemma 2.2]{LLWY}. In \cite{LLWY}, this lower bound is proved for $N\geq4$ by the change of variable (see also Nguyen \cite{H18}), and for $N=3$ by the change of variable $t=\sinh\rho$. We show in Section \ref{sec:TM} that the change of variable $t=\sinh\rho$ gives this lower bound in every dimension $N\geq3$, with an explicit weight.

In dimension two, Liang, Lu, Wang and Yang \cite[Theorem 1.7]{LLWY} also proved a stronger inequality under a weaker constraint. More precisely, they showed that for $0<\lambda<\frac14$ there is a constant $C>0$ such that $\int_{\mathbb{H}^2}\Phi_2(4\pi|u|^2)dV_{\mathbb{H}}\leq C$ for all $u\in W^{1,2}_0(\mathbb{H}^2)$ with
\[
\int_{\mathbb{H}^2}|\nabla_{\mathbb{H}}u|^{2}dV_{\mathbb{H}}-\frac14\int_{\mathbb{H}^2}|u|^{2}dV_{\mathbb{H}}-\lambda\int_{\mathbb{H}^2}(1-|x|^2)|u|^{2}dV_{\mathbb{H}}\leq1 .
\]
The remainder term in Theorem \ref{thm:pge3} gives a stronger inequality of this type in every dimension $N\geq3$.

\begin{theorem}\label{thm:HTMN}
Let $N\geq3$ and $0\leq\lambda<1$. Then there is a constant $C=C(N,\lambda)$ such that
\[
\int_{\mathbb{H}^N}\Phi_{N}\left(\beta_{N}|u|^{\frac{N}{N-1}}\right)dV_{\mathbb{H}}\leq C
\qquad\text{and}\qquad
\int_{B^N}e^{\beta_{N}|u|^{\frac{N}{N-1}}}dx\leq C
\]
for all $u\in W^{1,N}_{0}(\mathbb{H}^N)$ with
\begin{equation}\label{eq:constraintN}
\int_{\mathbb{H}^N}|\nabla_{\mathbb{H}}u|^{N}dV_{\mathbb{H}}-\left(\frac{N-1}{N}\right)^{N}\int_{\mathbb{H}^N}|u|^{N}dV_{\mathbb{H}}
-\lambda\left(\frac{N-1}{N}\right)^{N}\int_{\mathbb{H}^N}\frac{2}{e^{2\rho}-1}\,|u|^{N}dV_{\mathbb{H}}\leq 1 .
\end{equation}
\end{theorem}

For $\lambda=0$, Theorem \ref{thm:HTMN} is the result of Liang, Lu, Wang and Yang. For $0<\lambda<1$, the constraint \eqref{eq:constraintN} is weaker than the constraint in \eqref{eq:HTM}, since the weight $\frac{2}{e^{2\rho}-1}$ is positive. So Theorem \ref{thm:HTMN} is stronger than the result in \cite{LLWY} in every dimension $N\geq3$. The weight $\frac{2}{e^{2\rho}-1}$ behaves like $\rho^{-1}$ near the origin and it is exponentially small at infinity. In dimension three, the weight $R_{3}$ gives the following result, which is even stronger. Indeed, since $\frac{2}{e^{2\rho}-1}=\coth\rho-1$, we have
\[
\frac49R_3(\rho)=\frac{8}{27}\cdot\frac{2}{e^{2\rho}-1}+\frac{16}{27}\left[\frac{1}{\rho^{3}}-\frac{\coth\rho}{\sinh^{2}\rho}\right]\geq\frac{8}{27}\cdot\frac{2}{e^{2\rho}-1},
\]
since the bracket is nonnegative by the inequality $\sinh^3\rho\geq\rho^3\cosh\rho$ of Lazarevi\'c, see \eqref{eq:Laz} below. So the constraint \eqref{eq:constraint3} below is weaker than \eqref{eq:constraintN} with $N=3$. Moreover, as $\rho\to\infty$, the weight $\frac49R_3(\rho)$ behaves like $\frac{16}{27\rho^3}$, which is much larger than $\frac{8}{27}\cdot\frac{2}{e^{2\rho}-1}$.

\begin{theorem}\label{thm:HTM3}
Let $0\leq\lambda<1$ and let $R_{3}$ be the function \eqref{eq:beta} with $N=p=3$, that is,
\[
R_{3}(r)=\frac43\left[\frac{1}{r^{3}}-\frac{\coth r}{\sinh^{2}r}\right]+\frac23\left(\coth r-1\right).
\]
Then there is a constant $C=C(\lambda)$ such that
\[
\int_{\mathbb{H}^3}\Phi_{3}\left(\beta_{3}|u|^{3/2}\right)dV_{\mathbb{H}}\leq C
\qquad\text{and}\qquad
\int_{B^3}e^{\beta_{3}|u|^{3/2}}dx\leq C
\]
for all $u\in W^{1,3}_{0}(\mathbb{H}^3)$ with
\begin{equation}\label{eq:constraint3}
\int_{\mathbb{H}^3}|\nabla_{\mathbb{H}}u|^{3}dV_{\mathbb{H}}-\frac{8}{27}\int_{\mathbb{H}^3}|u|^{3}dV_{\mathbb{H}}
-\lambda\,\frac49\int_{\mathbb{H}^3}R_{3}(\rho)\,|u|^{3}dV_{\mathbb{H}}\leq 1 .
\end{equation}
\end{theorem}

The constraint \eqref{eq:constraint3} is weaker than the constraint in \eqref{eq:HTM}, since $R_{3}\geq0$.

\medskip
\subsection{Ideas of the proofs.}\label{sec:methods}
The three known proofs of \eqref{eq:PS} for $p=2$ and $N=3$ use tools which only work for the Laplacian. These tools are a Green's function, a Fourier transform, and an Euler--Lagrange equation whose solutions can be classified. We describe them briefly, since this helps to understand our method for $p\neq2$.

\emph{The proof of Benguria, Frank and Loss \cite{BFL08}.} They work in the half space model $\mathbb{R}^{3}_{+}=\{(x,y):x\in\mathbb{R}^{2},\,y>0\}$ of $\mathbb{H}^3$. In this model, \eqref{eq:PS} with $N=3$ and $p=2$ becomes the Hardy--Sobolev--Maz'ya inequality
\[
\int_{\mathbb{R}^{3}_{+}}|\nabla f|^{2}\,dx\,dy-\frac14\int_{\mathbb{R}^{3}_{+}}\frac{f^{2}}{y^{2}}\,dx\,dy
\geq S(3,2)^{2}\left(\int_{\mathbb{R}^{3}_{+}}|f|^{6}dx\,dy\right)^{1/3}.
\]
Maz'ya asked for the sharp constant in this inequality. The proof is by duality. The inequality is equivalent to a Hardy--Littlewood--Sobolev type inequality for the kernel of $\left(-\Delta-\frac{1}{4y^{2}}\right)^{-1/2}$. The change of function $f=\sqrt{y}\,g$ turns $-\Delta-\frac{1}{4y^{2}}$ into the Laplacian in $\mathbb{R}^4$ acting on functions which are radial in the last two variables. So the heat kernel of $-\Delta-\frac{1}{4y^{2}}$ can be computed, and the kernel of $\left(-\Delta-\frac{1}{4y^{2}}\right)^{-s}$, $0<s<\frac{3}{2}$, is bounded by the Riesz kernel of $(-\Delta)^{-s}$ in $\mathbb{R}^3$. Then the sharp Hardy--Littlewood--Sobolev inequality of Lieb \cite{Lieb83} gives the result.

\emph{The proof of Mancini and Sandeep \cite{MS08}.} They study the equation
\[
-\Delta_{\mathbb{H}}u-\lambda u=u^{q-1}\qquad\text{on }\mathbb{H}^N,
\]
with $\lambda\leq\frac{(N-1)^{2}}{4}$ and $2<q\leq\frac{2N}{N-2}$. They prove existence, uniqueness and nonexistence of positive solutions with finite energy. For $N=3$ and $q=6$, they prove that the equation has no positive solution with finite energy for every $\lambda\leq1$. Their proof uses a moving plane argument, which shows that positive solutions are radial, and a Pohozaev type identity for the radial equation. If the best constant in \eqref{PoinSob} with $N=3$ were smaller than $S(3,2)$, it would be attained by concentration compactness, and this would give a positive solution. So $C(3)=S(3,2)$, and the constant is not attained. For $N\geq4$, they prove that a positive solution of the critical equation exists if and only if $\frac{N(N-2)}{4}<\lambda\leq\frac{(N-1)^{2}}{4}$. This is the theorem of Hebey seen from the equation.

\emph{The proof of Lu and Yang \cite{LY19, LY221}.} They use the Helgason--Fourier transform on $\mathbb{H}^N$. With this transform they compute the Green's function of $-\Delta_{\mathbb{H}}-\frac{(N-1)^{2}}{4}$, and more generally of the Paneitz and GJMS operators, as explicit functions of the distance. They also prove a Hardy--Littlewood--Sobolev inequality on $\mathbb{H}^N$ for the kernel $\left(2\sinh\frac{\rho(x,y)}{2}\right)^{-\lambda}$, $0<\lambda<N$, with the sharp Euclidean constant \cite[Theorem 4.1]{LY19}. Then \eqref{eq:PS} follows by duality from a pointwise bound of the Green's function by this kernel. In this way they prove the sharp $k$-th order inequality on $\mathbb{H}^{2k+1}$, and the higher order Hardy--Sobolev--Maz'ya inequalities on half spaces. In \cite{LY221} they find precise formulas and best pointwise bounds for these Green's functions, and they show that for $N\geq2k+2$ the sharp constant is smaller than the Sobolev constant. For $k=1$ this gives the theorem of Hebey again. The same method also works on the complex, quaternionic and Cayley hyperbolic spaces \cite{LY22, FLY24}.

\emph{Our method.} For $p\neq2$ none of these tools is available. There is no Green's function, no Fourier transform for the $p$-Laplacian, and no classification of the solutions of the critical equation. The rearrangement method of Nguyen \cite{H18}, which was the only method for $p\neq2$, needs $p\geq\frac{2N}{N-1}$. In dimension three this means $p\geq3=N$, so it gives nothing for \eqref{eq:PS}. Our method uses the one part of the proofs for $p=2$ which does not depend on the Laplacian. This is the change of function. For $p=2$, the map $u\mapsto\left(\frac{\sinh r}{r}\right)^{\frac{N-1}{2}}u$ is a unitary map from $L^{2}(\sinh^{N-1}r\,dr)$ onto $L^{2}(r^{N-1}dr)$. This map has been used for the nonlinear Schr\"odinger equation on $\mathbb{H}^N$ by Christianson and Marzuola \cite{CM10} and by Banica and Duyckaerts \cite{BD15}. For $N=3$, it gives the identity
\[
\int_0^\infty\left(|u'|^{2}-u^{2}\right)\sinh^{2}r\,dr=\int_0^\infty|v'|^{2}r^{2}dr,\qquad v=u\,\frac{\sinh r}{r},
\]
which is behind \cite{BFL08} and \cite[Lemma 2.2]{MS08}. We use the $L^{p}$ version of this map,
\begin{equation}\label{eq:mapintro}
v(r)=u(r)\left(\frac{\sinh r}{r}\right)^{\frac{N-1}{p}} .
\end{equation}
With this map we have $u^{p}\sinh^{N-1}r=v^{p}r^{N-1}$. So the $L^{p}$ norm of $u$ on $\mathbb{H}^N$ is equal to the $L^{p}$ norm of $v$ on $\mathbb{R}^N$, and the $L^{q}$ norm of $v$ is larger than the $L^{q}$ norm of $u$ for $q>p$. For $p\neq2$ there is no identity like the one above. In place of the identity we use a one-dimensional convexity inequality, which is the tangent line inequality of Lemma \ref{lem:tangent}. The point of contact of the tangent line corresponds to the ground state $\sinh^{-\frac{N-1}{p}}\rho$ of $\mathbb{H}^N$. With this choice, the first order term of the tangent line cancels with the term from the integration by parts. What is left is the remainder term with the weight $R_{N}$. So everything is reduced to the sign of the explicit function $R_{N}$, which we prove by showing that $R_{N}$ is monotone. The monotonicity is also what we need in the Hardy--Littlewood inequality, when we pass from radial functions to all functions. The proof is elementary and short, and it works in every dimension and for every $2\leq p<N+2$, as long as $R_{N}$ is nonnegative. It gives for the first time the inequality on $\mathbb{H}^3$ for all $2<p<3$ and on $\mathbb{H}^4$ for $\frac{1+\sqrt{13}}{2}\leq p<\frac83$. It also gives a remainder which decays like $\rho^{-p}$ at infinity. Since the map keeps the $L^{p}$ norm, the same proof gives the sharp logarithmic Sobolev inequality. For $p\geq3$ we use the change of variable $t=\sinh\rho$ of \cite{LLWY} and the Hardy inequality. This gives Theorem \ref{thm:pge3}. At the endpoint $p=N$, it also gives a short proof, the same in every dimension $N\geq3$, of the key estimate in the proof of the Hardy--Trudinger--Moser inequality of Liang, Lu, Wang and Yang, and Theorem \ref{thm:HTMN}.

For Theorem \ref{thm:gap} we work on $\mathbb{H}^N$ itself. We write $u=g\sinh^{-\frac{N-1}{p}}\rho$. Then the $p$-Dirichlet integral minus the full spectral gap is equal to an integral of a positive weight times $g^{p}$, plus the integral of the difference between $|t|^{p}$ and its tangent line, see Proposition \ref{prop:ground}. For $1<p<2$ the function $|t|^{p}$ is not twice differentiable at $t=0$, and this difference is bounded from below by $\frac{p(p-1)}{2}|g'|^{2}\left(\frac{N-1}{p}\coth\rho\,g+|g'|\right)^{p-2}$, as in \cite{BKS26}. These two terms control the $L^{\frac{Np}{N-p}}$ norm of $u$ near the origin by the Sobolev inequality, and far from the origin by a bound on the growth of $g$.

\medskip
The paper is organized as follows.
\begin{itemize}
\item[Section 2:] We introduce the change of function and the tangent line inequality. Then we prove the key inequality with the weight $R_{N}$, and we study the sign and the monotonicity of $R_{N}$.
\item[Section 3:] We prove Theorem \ref{thm:main}.
\item[Section 4:] We use the change of variable $t=\sinh\rho$ to prove Theorem \ref{thm:pge3}. Then we prove Theorem \ref{thm:PS}.
\item[Section 5:] We use the ground state transform to prove Theorem \ref{thm:gap}.
\item[Section 6:] We prove the logarithmic Sobolev inequalities in Theorem \ref{thm:logPS} and Theorem \ref{thm:log}.
\item[Section 7:] We study the endpoint $p=N$ and prove Theorem \ref{thm:HTMN} and Theorem \ref{thm:HTM3}.
\item[Section 8:] We show that the inequalities \eqref{eq:PS} and \eqref{eq:logPS} fail for some $p<2$.
\end{itemize}

\medskip

\section{The change of function and the key inequality}\label{sec:key}

In this section we introduce the change of function, we prove the tangent line inequality and the key inequality, and we study the function $R_{N}$. To make the formulas shorter, in the rest of the paper we write
\[
\beta=\frac{N-1}{p}.
\]
So the constant in the Poincar\'e inequality \eqref{eq:gap} is $\beta^p$, and the remainder terms in Theorem \ref{thm:main} and Theorem \ref{thm:pge3} are $\beta^{p-1}\int R_N(\rho)|u|^pdV_{\mathbb{H}}$ and $\beta^{p}\int \frac{2}{e^{2\rho}-1}|u|^pdV_{\mathbb{H}}$.

\subsection{The change of function.}\label{sec:map}
We write $r$ for the geodesic distance to the origin. For a radial function $u=u(r)$ on $\mathbb{H}^N$, we define
\begin{equation}\label{eq:map}
v(r):=u(r)\,m(r),\qquad m(r):=\left(\frac{\sinh r}{r}\right)^{\beta},
\end{equation}
and we see $v$ as a radial function on $\mathbb{R}^N$. So the point at hyperbolic distance $r$ from the origin goes to the point at Euclidean distance $r$ from the origin, and only the value of the function is changed. The exponent $\beta$ is chosen so that
\begin{equation}\label{eq:weights}
m(r)^{p}\,r^{N-1}=\sinh^{N-1}r .
\end{equation}
Since $u=v/m$ and $\frac{m'}{m}=\beta\left(\coth r-\frac1r\right)$, we have
\[
u'=\frac{1}{m}\left(v'-\beta\left(\coth r-\frac1r\right)v\right).
\]
So, by \eqref{eq:weights},
\begin{equation}\label{eq:transfer}
|u'(r)|^{p}\sinh^{N-1}r=r^{N-1}\left|v'(r)-\beta\left(\coth r-\frac1r\right)v(r)\right|^{p},
\qquad
|u(r)|^{p}\sinh^{N-1}r=|v(r)|^{p}\,r^{N-1}.
\end{equation}
The second identity shows that the $L^p$ norm does not change. For the critical exponent, since $\beta \frac{Np}{N-p}>\beta p=N-1$ and $\sinh r\geq r$, we have
\begin{equation}\label{eq:norm}
\begin{aligned}
\int_{\mathbb{R}^N}|v|^{\frac{Np}{N-p}}dx
&=|\mathbb{S}^{N-1}|\int_{0}^{\infty}|u|^{\frac{Np}{N-p}}\left(\frac{\sinh r}{r}\right)^{\beta \frac{Np}{N-p}-(N-1)}\sinh^{N-1}r\,dr\\
&\geq \int_{\mathbb{H}^N}|u|^{\frac{Np}{N-p}}dV_{\mathbb{H}},
\end{aligned}
\end{equation}
and the inequality is strict if $u\not\equiv0$.

\subsection{A one-dimensional inequality.}\label{sec:1d}

\begin{lemma}\label{lem:tangent}
Let $p\geq2$, $\xi>0$ and $\eta_{0}\leq0$. Then for all $\eta\leq\xi$, we have
\begin{equation}\label{eq:tangent}
(\xi-\eta)^{p}-|\eta|^{p}\geq(\xi-\eta_{0})^{p}-|\eta_{0}|^{p}
-p\left[(\xi-\eta_{0})^{p-1}-|\eta_{0}|^{p-1}\right](\eta-\eta_{0}).
\end{equation}
\end{lemma}

The right side of \eqref{eq:tangent} is the tangent line at $\eta_{0}$ of the function $f(\eta)=(\xi-\eta)^{p}-|\eta|^{p}$. This function is not convex on $(-\infty,\xi]$, and the lemma is false for some $\eta_{0}>0$. So the condition $\eta_{0}\leq0$ is important. For $p=2$, the two sides of \eqref{eq:tangent} are equal.

\begin{proof}
Let $f(\eta)=(\xi-\eta)^{p}-|\eta|^{p}$ on $(-\infty,\xi]$, and let $\ell$ be the tangent line of $f$ at $\eta_{0}$. We have
\[
f''(\eta)=p(p-1)\left[(\xi-\eta)^{p-2}-|\eta|^{p-2}\right].
\]
Since $p\geq2$, $f$ is convex on $(-\infty,\xi/2]$ and concave on $[\xi/2,\xi]$. Since $\eta_{0}\leq0<\xi/2$, the convexity gives $f\geq\ell$ on $(-\infty,\xi/2]$. In particular, $f(\xi/2)\geq\ell(\xi/2)$. On $[\xi/2,\xi]$, the concave function $f$ is above the line segment between the points $(\xi/2,f(\xi/2))$ and $(\xi,f(\xi))$. So $f\geq \ell$ on $[\xi/2,\xi]$ if $\ell(\xi)\leq f(\xi)=-\xi^{p}$. Therefore it is enough to prove
\[
\ell(\xi)\leq-\xi^{p}.
\]
Let $a=\xi-\eta_{0}$ and $b=-\eta_{0}$. Then $0\leq b<a$ and $\xi=a-b$. We have
\[
\ell(\xi)=a^{p}-b^{p}-p\left(a^{p-1}-b^{p-1}\right)a=-(p-1)a^{p}-b^{p}+p\,b^{p-1}a .
\]
We divide by $a^{p}$ and write $t=b/a\in[0,1)$. Then the inequality $\ell(\xi)\leq-(a-b)^{p}$ becomes
\begin{equation}\label{eq:Phi}
\Psi(p,t):=(p-1)+t^{p}-p\,t^{p-1}-(1-t)^{p}\geq0 .
\end{equation}
We have $\Psi(2,t)=1+t^{2}-2t-(1-t)^{2}=0$, and for $p\geq2$ and $0\leq t<1$,
\[
\frac{\partial\Psi}{\partial p}(p,t)=\left(1-t^{p-1}\right)+t^{p-1}(t-p)\ln t-(1-t)^{p}\ln(1-t)\geq0,
\]
since each of the three terms is nonnegative. So $\Psi(p,t)\geq\Psi(2,t)=0$. This proves \eqref{eq:Phi} and the lemma.
\end{proof}

\subsection{The key inequality.}\label{sec:keyineq}

\begin{proposition}\label{prop:key}
Let $N\geq3$ and $2\leq p<N+2$. Let $u\geq0$ be a radial, nonincreasing and Lipschitz function with compact support, and let $v$ be as in \eqref{eq:map}. Then
\begin{equation}\label{eq:key}
\int_{0}^{\infty}\left(|u'|^{p}-\beta^{p}u^{p}\right)\sinh^{N-1}r\,dr
\geq\int_{0}^{\infty}|v'|^{p}\,r^{N-1}dr
+\beta^{p-1}\int_{0}^{\infty}R_{N}(r)\,v^{p}\,r^{N-1}dr ,
\end{equation}
where $R_{N}$ is the function in \eqref{eq:beta}.
\end{proposition}

\begin{proof}
We apply Lemma \ref{lem:tangent} at each $r>0$ with
\begin{equation}\label{eq:choice}
\xi:=\beta\left(\coth r-\frac1r\right)v(r),\qquad
\eta:=v'(r),\qquad
\eta_{0}:=-\beta\,\frac{v(r)}{r}\,.
\end{equation}
The conditions of the lemma hold. Indeed, $\eta_{0}\leq0$ since $v\geq0$, and $\eta\leq\xi$ since
\[
v'-\beta\left(\coth r-\frac1r\right)v=m\,u'\leq0 .
\]
At the points where $v=0$, we have $\xi=0$, and \eqref{eq:tangent} becomes $|\eta|^{p}\leq|\eta|^{p}$. Also,
\[
\xi-\eta_{0}=\beta\,(\coth r)\,v,\qquad
|\eta_{0}|=\beta\,\frac{v}{r},\qquad
\eta-\eta_{0}=v'+\beta\,\frac{v}{r}.
\]
We multiply \eqref{eq:tangent} by $r^{N-1}$ and use the first identity in \eqref{eq:transfer}. We get
\[
\begin{aligned}
|u'|^{p}\sinh^{N-1}r-|v'|^{p}r^{N-1}\geq
r^{N-1}\Bigg\{&\beta^{p}\left(\coth^{p}r-\frac{1}{r^{p}}\right)v^{p}\\
&-p\,\beta^{p-1}\left(\coth^{p-1}r-\frac{1}{r^{p-1}}\right)v^{p-1}
\left(v'+\beta\,\frac{v}{r}\right)\Bigg\}.
\end{aligned}
\]
We subtract $\beta^{p}u^{p}\sinh^{N-1}r=\beta^{p}v^{p}r^{N-1}$ from both sides. In the last term we use $p\beta=N-1$ and $v^{p-1}v'=\frac1p(v^{p})'$. Then we get
\begin{equation}\label{eq:pointwise}
\begin{aligned}
\left(|u'|^{p}-\beta^{p}u^{p}\right)\sinh^{N-1}r-|v'|^{p}r^{N-1}
\geq\ &\beta^{p}\left(\coth^{p}r-\frac{1}{r^{p}}-1\right)v^{p}\,r^{N-1}\\
&-(N-1)\,\beta^{p-1}\left(\coth^{p-1}r-\frac{1}{r^{p-1}}\right)v^{p}\,r^{N-2}\\
&-\beta^{p-1}\left(\coth^{p-1}r-\frac{1}{r^{p-1}}\right)r^{N-1}\,\left(v^{p}\right)' .
\end{aligned}
\end{equation}
Now we integrate \eqref{eq:pointwise} over $(0,\infty)$, and we integrate the last term by parts. The boundary terms are zero. This is because $u$ has compact support, and as $r\to0^{+}$,
\[
\coth^{p-1}r-\frac{1}{r^{p-1}}=\frac{p-1}{3}r^{3-p}+O(r^{5-p}),\qquad\text{so}\qquad \left(\coth^{p-1}r-\frac{1}{r^{p-1}}\right)r^{N-1}=O(r^{N+2-p})\to0,
\]
since $p<N+2$. For the same reason, all the integrals below are finite.
Since $(\coth r)'=-\frac{1}{\sinh^{2}r}$, we have
\[
\left[\left(\coth^{p-1}r-\frac{1}{r^{p-1}}\right)r^{N-1}\right]'
=(p-1)\left(\frac{1}{r^{p}}-\frac{\coth^{p-2}r}{\sinh^{2}r}\right)r^{N-1}
+(N-1)\left(\coth^{p-1}r-\frac{1}{r^{p-1}}\right)r^{N-2} .
\]
The second term on the right cancels the second line of \eqref{eq:pointwise}. This is the reason for our choice of $\eta_0$ in \eqref{eq:choice}. Therefore,
\[
\int_{0}^{\infty}\left(|u'|^{p}-\beta^{p}u^{p}\right)\sinh^{N-1}r\,dr
-\int_{0}^{\infty}|v'|^{p}r^{N-1}dr\geq\int_{0}^{\infty}\Lambda(r)\,v^{p}\,r^{N-1}dr,
\]
where
\[
\Lambda(r):=\beta^{p}\left(\coth^{p}r-\frac{1}{r^{p}}-1\right)
+(p-1)\beta^{p-1}\left(\frac{1}{r^{p}}-\frac{\coth^{p-2}r}{\sinh^{2}r}\right).
\]
Finally, since
\[
\coth^{p}r-1=\coth^{p-2}r\left(\coth^{2}r-1\right)+\left(\coth^{p-2}r-1\right)
=\frac{\coth^{p-2}r}{\sinh^{2}r}+\left(\coth^{p-2}r-1\right),
\]
we have
\[
\begin{aligned}
\beta^{1-p}\Lambda(r)
&=\beta\left(\coth^{p}r-\frac{1}{r^{p}}-1\right)
+(p-1)\left(\frac{1}{r^{p}}-\frac{\coth^{p-2}r}{\sinh^{2}r}\right)\\
&=(p-1-\beta)\left(\frac{1}{r^{p}}-\frac{\coth^{p-2}r}{\sinh^{2}r}\right)
+\beta\left(\coth^{p-2}r-1\right)=R_{N}(r).
\end{aligned}
\]
This proves \eqref{eq:key}.
\end{proof}

When $N=3$ and $p=2$, both terms of $R_{3}$ are zero, so $R_{3}\equiv0$ and \eqref{eq:key} is an identity. In fact, Lemma \ref{lem:tangent} is an identity for $p=2$, and \eqref{eq:key} becomes
\[
\int_0^\infty\left(|u'|^{2}-u^{2}\right)\sinh^{2}r\,dr=\int_0^\infty|v'|^{2}r^{2}dr,\qquad v=u\,\frac{\sinh r}{r},
\]
which is behind the inequality of Benguria, Frank and Loss.

\subsection{The function $R_{N}$.}\label{sec:B}
In this subsection we write
\[
S=S(r)=\frac{\sinh r}{r},\qquad C=C(r)=\cosh r .
\]
We will use the inequality of Lazarevi\'c \cite[3.6.9]{Mit70},
\begin{equation}\label{eq:Laz}
\sinh^{3}r\geq r^{3}\cosh r,\qquad\text{that is,}\qquad S^{3}\geq C,\qquad r>0 .
\end{equation}
For the convenience of the reader we give a short proof. Let $f(r)=3\ln\frac{\sinh r}{r}-\ln\cosh r$. Then $f(0^+)=0$ and
\[
f'(r)=3\coth r-\frac3r-\tanh r .
\]
With $s=2r$, the inequality $f'(r)\geq0$ is equivalent to $s\cosh s+2s-3\sinh s\geq0$, and the left side is equal to
\[
\sum_{k\geq2}\frac{2k-2}{(2k+1)!}\,s^{2k+1}\geq 0 .
\]
So $f\geq0$, which is \eqref{eq:Laz}.

\begin{lemma}\label{lem:mono}
Let $N\geq3$ and $p\geq2$. Assume that $p-1\geq\beta$, that is, $p^{2}-p\geq N-1$ (this is always true for $N=3$). Assume also that one of the following conditions holds:
\begin{itemize}
\item[(a)] $p\leq3$,
\item[(b)] $3<p\leq5$ and $2p\leq N+1$.
\end{itemize}
Then $R_{N}$ is nonincreasing on $(0,\infty)$ and $R_{N}(r)\to0$ as $r\to\infty$. In particular, $R_{N}\geq0$. If $(N,p)\neq(3,2)$, then $R_{N}$ is strictly decreasing and strictly positive, and \eqref{eq:Basym} holds.
\end{lemma}

The cases (i)--(iii) of Theorem \ref{thm:main} are exactly the cases of this lemma with $p<N$. Indeed, for $N=3$ the condition $p-1\geq\beta$ holds for all $p\geq2$, and (a) with $p<3$ is (i). For $N=4$ we have $2p\leq N+1$ only if $p\leq\frac52$, so (b) is empty and (a) gives (ii). For $N\geq5$ and $p\leq3$, we have $2p\leq6\leq N+1$. So (a) and (b) together give $p\leq\min\{5,\frac{N+1}{2}\}$, which is (iii).

\begin{proof}
\emph{Step 1: a formula for $R_{N}'$.} Let $c_{1}=p-1-\beta\geq0$, $w=r\coth r$ and $\sigma=\frac{r^{2}}{\sinh^{2}r}$. Then
\[
w^{2}-\sigma=r^{2},\qquad w'=\frac{w-\sigma}{r},\qquad \sigma'=\frac{2\sigma(1-w)}{r}.
\]
With these functions we can write
\[
R_{N}=c_{1}\,r^{-p}\left(1-w^{p-2}\sigma\right)+\beta\left(w^{p-2}r^{2-p}-1\right).
\]
By a direct computation, we get
\begin{equation}\label{eq:Bprime}
R_{N}'(r)=\frac{1}{r^{p+1}}\left[-p\,c_{1}+(p-1)(p-2)\,w^{p-3}\sigma^{2}+(2p-N-1)\,w^{p-1}\sigma\right].
\end{equation}
Here we used $w^{2}-\sigma=r^{2}$, $c_{1}+\beta=p-1$ and $p\beta=N-1$. Since $w=C/S$ and $\sigma=S^{-2}$, we have
\[
w^{p-3}\sigma^{2}=C^{p-3}S^{-(p+1)},\qquad w^{p-1}\sigma=C^{p-1}S^{-(p+1)}.
\]
Using $C^{2}=1+\sinh^{2}r$ and $(p-1)(p-2)+(2p-N-1)=p\,c_{1}$, we get
\begin{equation}\label{eq:Bprime2}
r^{p+1}S^{p+1}C^{3-p}\,R_{N}'(r)
=-p\,c_{1}\left(S^{p+1}C^{3-p}-1\right)+(2p-N-1)\sinh^{2}r .
\end{equation}

\emph{Step 2: the case $2p\leq N+1$ and $p\leq5$.} In this case the last term in \eqref{eq:Bprime2} is nonpositive. We claim that the first term is also nonpositive, that is, $S^{p+1}C^{3-p}\geq1$. For $p\leq3$ this is clear, since $S\geq1$, $C\geq1$ and $3-p\geq0$. For $3<p\leq5$, we use \eqref{eq:Laz} in the form $S^{p+1}\geq C^{(p+1)/3}$. Since $\frac{p+1}{3}\geq p-3$ when $p\leq5$, we get $S^{p+1}\geq C^{p-3}$, which is the claim. So $R_{N}'\leq0$. This step covers the case (b), the case (a) for $N\geq5$ (since then $p\leq3\leq\frac{N+1}{2}$), and the case $N=4$, $p\leq\frac52$.

\emph{Step 3: the case $N=4$ and $\frac52<p\leq3$.} This is the rest of the case (a). Here $c_{1}=p-1-\frac3p>0$ and $2p-5\in(0,1]$. We claim that
\begin{equation}\label{eq:N4}
S^{p+1}C^{3-p}\geq1+\frac12\sinh^{2}r\qquad\text{for } \frac52\leq p\leq3 .
\end{equation}
If \eqref{eq:N4} holds, then by \eqref{eq:Bprime2},
\[
r^{p+1}S^{p+1}C^{3-p}R_{4}'\leq\left(-\frac{p\,c_{1}}{2}+2p-5\right)\sinh^{2}r
=-\frac{p^{2}-5p+7}{2}\sinh^{2}r<0,
\]
since $p^{2}-5p+7>0$ for all $p$. To prove \eqref{eq:N4}, we note that its logarithm
\[
(p+1)\ln S+(3-p)\ln C-\ln\left(1+\frac12\sinh^{2}r\right)
\]
is a linear function of $p$. So it is enough to check the two endpoints $p=3$ and $p=\frac52$. For $p=3$ we need
\begin{equation}\label{eq:p3}
S^{4}\geq1+\frac12\sinh^{2}r .
\end{equation}
With $s=\sinh^{2}r$, this is $s^{2}-\frac{r^{4}}{2}s-r^{4}\geq0$, and it holds if
\[
s\geq\frac{r^{4}}{4}+r^{2}\sqrt{1+\frac{r^{4}}{16}}.
\]
Since $\sqrt{1+t}\leq1+\frac t2$, the right side is at most $r^{2}+\frac{r^{4}}{4}+\frac{r^{6}}{32}$. The Taylor series of $\sinh^{2}r=\frac12(\cosh2r-1)$ has positive coefficients, so
\[
s=\sinh^{2}r\geq r^{2}+\frac{r^{4}}{3}+\frac{2r^{6}}{45}\geq r^{2}+\frac{r^{4}}{4}+\frac{r^{6}}{32}.
\]
This proves \eqref{eq:p3}. For $p=\frac52$ we need $S^{7/2}C^{1/2}\geq1+\frac12\sinh^{2}r$. This follows from \eqref{eq:p3}, since
\[
S^{7/2}C^{1/2}=S^{4}\left(\frac{C}{S}\right)^{1/2}=S^{4}\,(r\coth r)^{1/2}\geq S^{4}.
\]

\emph{Step 4: the case $N=3$.} Here $c_{1}=\frac{(p-2)(p+1)}{p}$ and $2p-N-1=2(p-2)$. So by \eqref{eq:Bprime2}, $R_{3}'\leq0$ if and only if
\begin{equation}\label{eq:star}
S^{p+1}C^{3-p}\geq1+\frac{2\sinh^{2}r}{p+1} .
\end{equation}
For $p=2$ there is nothing to prove, since $R_{3}\equiv0$. Fix $r>0$ and let
\[
g(p)=(p+1)\ln S+(3-p)\ln C-\ln\left(1+\frac{2\sinh^{2}r}{p+1}\right),\qquad 2\leq p\leq3 .
\]
The first two terms are linear in $p$. With $h(p)=\frac{2\sinh^{2}r}{p+1}$, we have
\[
\frac{d^{2}}{dp^{2}}\ln(1+h)=\frac{h''(1+h)-h'^{2}}{(1+h)^{2}}
=\frac{1}{(1+h)^{2}}\left[\frac{4\sinh^{2}r}{(p+1)^{3}}+\frac{4\sinh^{4}r}{(p+1)^{4}}\right]>0 .
\]
So $g$ is strictly concave in $p$, and $g\geq0$ on $[2,3]$ follows from $g(2)\geq0$ and $g(3)\geq0$. For $p=3$, \eqref{eq:star} is $S^{4}\geq1+\frac12\sinh^{2}r$, which is \eqref{eq:p3}. For $p=2$ we need $S^{3}C\geq1+\frac23\sinh^{2}r$. By \eqref{eq:Laz}, we have
\[
S^{3}C\geq C^{2}=1+\sinh^{2}r>1+\frac23\sinh^{2}r.
\]
Hence $R_{3}'\leq0$. For $2<p<3$ the inequality is strict, since $g(2)>0$ and $g$ is strictly concave.

\emph{Step 5: limits and positivity.} As $r\to\infty$, the terms $\coth^{p-2}r-1$ and $\sinh^{-2}r$ in \eqref{eq:beta} are $O(e^{-2r})$. So $R_{N}(r)\to0$ and $r^{p}R_{N}(r)\to c_{1}$. As $r\to0^{+}$, we have
\[
\frac{\coth^{p-2}r}{\sinh^{2}r}=\frac{1}{r^{p}}\left(1+\frac{p-3}{3}r^{2}+O(r^{4})\right),\qquad
\coth^{p-2}r=r^{2-p}\left(1+O(r^{2})\right).
\]
Since $(N,p)\neq(3,2)$, we have $p>2$, so $r^{p-2}\to0$. Hence $r^{p-2}R_{N}(r)\to\frac{c_{1}(3-p)}{3}+\beta$. Using $c_{1}=p-1-\beta$ and $p\beta=N-1$, we get
\[
\frac{c_{1}(3-p)}{3}+\beta=\frac{(p-1)(3-p)+N-1}{3}=\frac{N-(p-2)^{2}}{3},
\]
which is \eqref{eq:Basym}. This limit is nonnegative in all cases of the lemma. Indeed, for $p\leq3$ we have $(p-2)^{2}\leq1<N$. In the case (b) we have $N\geq2p-1$, so
\[
N-(p-2)^{2}\geq2p-1-(p-2)^{2}=-p^{2}+6p-5=(p-1)(5-p)\geq0.
\]
Finally, $R_{N}\geq0$ since $R_{N}$ is nonincreasing and tends to $0$. If $R_{N}$ is strictly decreasing, then $R_{N}>0$. This is the case in Steps 2 and 3 if $c_{1}>0$ or $2p<N+1$. The conditions $c_{1}=0$ and $2p=N+1$ cannot hold at the same time, since $p^{2}-p=N-1=2p-2$ gives $p=2$ and $N=3$. It is also the case in Step 4 for $p>2$.
\end{proof}

When $p\leq3$ and $p-1\geq\beta$, the inequality $R_{N}\geq0$ follows directly from \eqref{eq:Laz}. Indeed, the inequality $\frac{1}{r^{p}}\geq\frac{\coth^{p-2}r}{\sinh^{2}r}$ is the same as $\sinh r\geq r\cosh^{(p-2)/p}r$, which follows from $S\geq C^{1/3}$ since $\frac{p-2}{p}\leq\frac13$. Also, $\coth r\geq1$. The monotonicity of $R_N$ is the extra information that we need in the next section. For the inequalities without remainder, we only need the sign of $R_{N}$, and the sign can be controlled in a larger range. By the two limits in \eqref{eq:Basym}, the conditions $p-1\geq\beta$ and $(p-2)^{2}\leq N$ are necessary for $R_{N}\geq0$. The next lemma gives a sufficient condition of the same form.

\begin{lemma}\label{lem:pos}
Let $N\geq3$ and $p\geq5$, and assume that $(p-1)(p-2)\leq N-1\leq p(p-1)$, that is,
\[
\frac{1+\sqrt{4N-3}}{2}\leq p\leq\frac{3+\sqrt{4N-3}}{2}.
\]
Then $R_{N}\geq0$ on $(0,\infty)$.
\end{lemma}

\begin{proof}
Let $y=\frac1r$ and $x=\coth r$. Then $x>1$, $x>y$ and $\frac{1}{\sinh^{2}r}=x^{2}-1$. With $c_{1}=p-1-\beta$, we have
\[
R_{N}=c_{1}\left(y^{p}-x^{p}+x^{p-2}\right)+\beta\left(x^{p-2}-1\right)=(p-1)P+\beta Q,
\]
where
\[
P:=y^{p}-x^{p}+x^{p-2},\qquad Q:=x^{p}-y^{p}-1 .
\]
The condition of the lemma is $\frac{(p-1)(p-2)}{p}\leq\beta\leq p-1$.

\emph{The case $Q\leq0$.} Then $\beta Q\geq(p-1)Q$, and so
\[
R_{N}\geq(p-1)(P+Q)=(p-1)\left(x^{p-2}-1\right)\geq0.
\]

\emph{The case $Q>0$.} Then $\beta Q\geq\frac{(p-1)(p-2)}{p}Q$, and so
\[
R_{N}\geq\frac{p-1}{p}\left[pP+(p-2)Q\right]=\frac{2(p-1)}{p}\left[y^{p}-G(x)\right],
\qquad G(x):=x^{p}-\frac p2x^{p-2}+\frac{p-2}{2}.
\]
So it is enough to show that $G(x)\leq y^{p}$. We have $G(1)=0$ and
\[
G'(x)=p\,x^{p-3}\left(x^{2}-\frac{p-2}{2}\right).
\]
If $x^{2}\leq\frac{p-2}{2}$, then $G$ is nonincreasing on $[1,x]$, so $G(x)\leq G(1)=0<y^{p}$. If $x^{2}>\frac{p-2}{2}$, let $z=x^{-2}<\frac{2}{p-2}\leq\frac23$. The inequality \eqref{eq:Laz} says that $\frac1r\geq\frac{\cosh^{1/3}r}{\sinh r}$. Since $\cosh r=x\sinh r$, this is $y\geq\left(x(x^{2}-1)\right)^{1/3}$, so that $y^{p}\geq x^{p}(1-z)^{p/3}$. By the Bernoulli inequality, $(1-z)^{p/3}\geq1-\frac p3z$, since $\frac p3\geq1$. Hence
\[
\begin{aligned}
y^{p}-G(x)&\geq x^{p}\left[1-\frac p3z-1+\frac p2z-\frac{p-2}{2}z^{p/2}\right]
=x^{p}z\left[\frac p6-\frac{p-2}{2}\,z^{\frac{p-2}{2}}\right]\\
&\geq x^{p}z\left[\frac p6-\frac{p-2}{2}\left(\frac{2}{p-2}\right)^{\frac{p-2}{2}}\right].
\end{aligned}
\]
With $k=\frac{p-2}{2}\geq\frac32$, the last bracket is $\frac13\left[k+1-3k^{1-k}\right]$. The function $k\mapsto3k^{1-k}$ is decreasing for $k\geq1$ and it is equal to $\sqrt6<\frac52$ at $k=\frac32$, while $k+1$ is increasing. So the bracket is positive, and $G(x)<y^{p}$.
\end{proof}

Lemma \ref{lem:pos} covers the interval $\max\left\{5,\frac{1+\sqrt{4N-3}}{2}\right\}\leq p\leq\frac{3+\sqrt{4N-3}}{2}$. This interval is not empty for $N\geq13$, and it has length one for $N\geq21$. It is inside the range $p\geq\frac{2N}{N-1}$ of the comparison of Nguyen. By Lemma \ref{lem:mono}, Lemma \ref{lem:pos} and the paragraph before Lemma \ref{lem:pos}, we know that $R_{N}\geq0$ under the condition $p-1\geq\beta$ when $p\leq3$, or $p\leq\min\{5,\frac{N+1}{2}\}$, or $5\leq p\leq\frac{3+\sqrt{4N-3}}{2}$.

\medskip
\section{The remainder with the weight $R_{N}$ -- Proof of Theorem \ref{thm:main}}\label{sec:proof}

In this section we prove Theorem \ref{thm:main}. So $(N,p)$ satisfies (i), (ii) or (iii) of Theorem \ref{thm:main}. Then $p\geq2$, $p<N$, and Lemma \ref{lem:mono} can be applied. We first prove the inequality for radial functions, then for all functions. Then we prove that the constants are sharp and that $S(N,p)^p$ is not attained.

\medskip
\emph{Radial functions.} Let $u\geq0$ be a radial, nonincreasing and Lipschitz function with compact support. Let $v=u\,m$ be as in \eqref{eq:map}. Then $v$ is a Lipschitz function with compact support on $\mathbb{R}^N$. By polar coordinates, Proposition \ref{prop:key}, the Sobolev inequality \eqref{eq:sobolev} for $v$ and \eqref{eq:norm}, we have
\begin{equation}\label{eq:chain}
\begin{aligned}
\int_{\mathbb{H}^N}|\nabla_{\mathbb{H}}u|^{p}dV_{\mathbb{H}}-\beta^{p}\int_{\mathbb{H}^N}u^{p}dV_{\mathbb{H}}
&=|\mathbb{S}^{N-1}|\int_{0}^{\infty}\left(|u'|^{p}-\beta^{p}u^{p}\right)\sinh^{N-1}r\,dr\\
&\geq |\mathbb{S}^{N-1}|\left[\int_{0}^{\infty}|v'|^{p}r^{N-1}dr
+\beta^{p-1}\int_{0}^{\infty}R_{N}(r)\,v^{p}r^{N-1}dr\right]\\
&=\int_{\mathbb{R}^N}|\nabla v|^{p}dx+\beta^{p-1}\int_{\mathbb{H}^N}R_{N}(\rho)\,u^{p}dV_{\mathbb{H}}\\
&\geq S(N,p)^{p}\left(\int_{\mathbb{R}^N}|v|^{\frac{Np}{N-p}}dx\right)^{\frac{N-p}{N}}+\beta^{p-1}\int_{\mathbb{H}^N}R_{N}(\rho)\,u^{p}dV_{\mathbb{H}}\\
&\geq S(N,p)^{p}\left(\int_{\mathbb{H}^N}|u|^{\frac{Np}{N-p}}dV_{\mathbb{H}}\right)^{\frac{N-p}{N}}+\beta^{p-1}\int_{\mathbb{H}^N}R_{N}(\rho)\,u^{p}dV_{\mathbb{H}} .
\end{aligned}
\end{equation}
In the third line we used $v^{p}r^{N-1}=u^{p}\sinh^{N-1}r$, see \eqref{eq:transfer}. This is \eqref{eq:main} for such $u$. Note that the remainder is finite. Indeed, $R_{N}$ is bounded away from the origin, and near the origin $R_{N}(\rho)u^{p}$ behaves like $c\,\rho^{2-p}u(0)^{p}$, which is integrable with respect to $dV_{\mathbb{H}}=\sinh^{N-1}\rho\,d\rho\,d\sigma$.

\medskip
\emph{All functions.} Let $u\in C_{c}^{\infty}(\mathbb{H}^N)$ and let $u^{\sharp}$ be the symmetric decreasing rearrangement of $|u|$ with respect to the hyperbolic volume, with center at the origin. Then $u^{\sharp}$ is radial, nonincreasing and Lipschitz with compact support, and it has the same distribution function as $|u|$. We use three facts. The first one is the P\'olya--Szeg\"o principle on $\mathbb{H}^N$ \cite{Bae94}:
\[
\int_{\mathbb{H}^N}|\nabla_{\mathbb{H}}u^{\sharp}|^{p}dV_{\mathbb{H}}\leq\int_{\mathbb{H}^N}|\nabla_{\mathbb{H}}u|^{p}dV_{\mathbb{H}}.
\]
The second one is that the rearrangement keeps the $L^{q}$ norms:
\[
\int_{\mathbb{H}^N}(u^{\sharp})^{q}dV_{\mathbb{H}}=\int_{\mathbb{H}^N}|u|^{q}dV_{\mathbb{H}}\qquad\text{for all } q\geq1.
\]
The third one comes from Lemma \ref{lem:mono}. The weight $R_{N}(\rho)$ is a nonnegative, radial and nonincreasing function, so it is equal to its own symmetric decreasing rearrangement. So, by the Hardy--Littlewood inequality \cite[Theorem 3.4]{LL01}, whose proof by the layer cake formula works on $(\mathbb{H}^N,dV_{\mathbb{H}})$ with rearrangements about the origin,
\[
\int_{\mathbb{H}^N}R_{N}(\rho)\,|u|^{p}dV_{\mathbb{H}}\leq\int_{\mathbb{H}^N}R_{N}(\rho)\,(u^{\sharp})^{p}dV_{\mathbb{H}} .
\]
Now we apply \eqref{eq:chain} to $u^{\sharp}$ and we use the three facts. We get
\[
\begin{aligned}
\int_{\mathbb{H}^N}|\nabla_{\mathbb{H}}u|^{p}dV_{\mathbb{H}}-\beta^{p}\int_{\mathbb{H}^N}|u|^{p}dV_{\mathbb{H}}
&\geq\int_{\mathbb{H}^N}|\nabla_{\mathbb{H}}u^{\sharp}|^{p}dV_{\mathbb{H}}-\beta^{p}\int_{\mathbb{H}^N}(u^{\sharp})^{p}dV_{\mathbb{H}}\\
&\geq S(N,p)^{p}\left(\int_{\mathbb{H}^N}(u^{\sharp})^{\frac{Np}{N-p}}dV_{\mathbb{H}}\right)^{\frac{N-p}{N}}+\beta^{p-1}\int_{\mathbb{H}^N}R_{N}(\rho)\,(u^{\sharp})^{p}dV_{\mathbb{H}}\\
&\geq S(N,p)^{p}\left(\int_{\mathbb{H}^N}|u|^{\frac{Np}{N-p}}dV_{\mathbb{H}}\right)^{\frac{N-p}{N}}+\beta^{p-1}\int_{\mathbb{H}^N}R_{N}(\rho)\,|u|^{p}dV_{\mathbb{H}} .
\end{aligned}
\]
This is \eqref{eq:main} for $u\in C_{c}^{\infty}(\mathbb{H}^N)$. For $u\in D^{1,p}(\mathbb{H}^N)$, we take $u_{n}\in C_{c}^{\infty}(\mathbb{H}^N)$ such that $u_{n}\to u$ in $D^{1,p}(\mathbb{H}^N)$ and almost everywhere. The left side of \eqref{eq:main} and the $L^{\frac{Np}{N-p}}$ norm are continuous in $D^{1,p}(\mathbb{H}^N)$, by \eqref{eq:gap} and by \eqref{eq:PS}, which we already proved for functions in $C_{c}^{\infty}(\mathbb{H}^N)$. By Fatou's lemma,
\[
\int_{\mathbb{H}^N}R_{N}(\rho)|u|^{p}dV_{\mathbb{H}}\leq\liminf_{n\to\infty}\int_{\mathbb{H}^N}R_{N}(\rho)|u_{n}|^{p}dV_{\mathbb{H}}.
\]
So \eqref{eq:main} holds for $u$.

\medskip
\emph{The constant $\beta^{p}$ is sharp.} Assume that \eqref{eq:PS} holds with some $c>\beta^{p}$ in place of $\beta^{p}$. Then
\[
\int_{\mathbb{H}^N}|\nabla_{\mathbb{H}}u|^{p}dV_{\mathbb{H}}\geq c\int_{\mathbb{H}^N}|u|^{p}dV_{\mathbb{H}}
\]
for all $u\in C_{c}^{\infty}(\mathbb{H}^N)$. This is not possible, since the constant in \eqref{eq:gap} is sharp.

\medskip
\emph{The constant $S(N,p)^{p}$ is sharp.} Let $w\in C_{c}^{\infty}(\mathbb{R}^N)$ with support in $\{|y|\leq R\}$. For small $\varepsilon>0$, we define $u_{\varepsilon}$ on $\mathbb{H}^N$ in geodesic normal coordinates at the origin by
\[
u_{\varepsilon}(\exp_{0}(\varepsilon y))=\varepsilon^{-\frac{N-p}{p}}w(y).
\]
In these coordinates the metric is $d\rho^{2}+\sinh^{2}\rho\,d\sigma^{2}$. On the support of $u_{\varepsilon}$ we have $\rho\leq\varepsilon R$, so $\frac{\sinh\rho}{\rho}=1+O(\varepsilon^{2})$. Hence
\[
\begin{gathered}
\int_{\mathbb{H}^N}|\nabla_{\mathbb{H}}u_{\varepsilon}|^{p}dV_{\mathbb{H}}=\left(1+O(\varepsilon^{2})\right)\int_{\mathbb{R}^N}|\nabla w|^{p}dy,\qquad
\int_{\mathbb{H}^N}|u_{\varepsilon}|^{p}dV_{\mathbb{H}}=O(\varepsilon^{p}),\\
\int_{\mathbb{H}^N}|u_{\varepsilon}|^{\frac{Np}{N-p}}dV_{\mathbb{H}}=\left(1+O(\varepsilon^{2})\right)\int_{\mathbb{R}^N}|w|^{\frac{Np}{N-p}}dy .
\end{gathered}
\]
So, as $\varepsilon\to0$, the quotient of the left side of \eqref{eq:PS} by $\left(\int_{\mathbb{H}^N}|u_{\varepsilon}|^{\frac{Np}{N-p}}dV_{\mathbb{H}}\right)^{\frac{N-p}{N}}$ tends to
\[
\frac{\int_{\mathbb{R}^N}|\nabla w|^{p}dy}{\left(\int_{\mathbb{R}^N}|w|^{\frac{Np}{N-p}}dy\right)^{\frac{N-p}{N}}}.
\]
Since $C_{c}^{\infty}(\mathbb{R}^N)$ is dense in $D^{1,p}(\mathbb{R}^N)$, this quotient can be as close to $S(N,p)^{p}$ as we want. So the constant $S(N,p)^{p}$ in \eqref{eq:PS}, and so in \eqref{eq:main}, cannot be larger.

\medskip
\emph{The constant $S(N,p)^{p}$ is not attained.} Assume that $u\in D^{1,p}(\mathbb{H}^N)$, $u\not\equiv0$, gives equality in \eqref{eq:PS}. Then $u^{\sharp}$ also gives equality, since the left side of \eqref{eq:PS} does not increase and the right side does not change. We take radial, nonincreasing and Lipschitz functions $u_{n}\geq0$ with compact support such that $u_{n}\to u^{\sharp}$ in $D^{1,p}(\mathbb{H}^N)$ and almost everywhere. Let $v_{n}=u_{n}m$ and $v=u^{\sharp}m$. By \eqref{eq:chain}, $\int_{\mathbb{R}^N}|\nabla v_{n}|^{p}dx$ is bounded by the left side of \eqref{eq:PS} for $u_{n}$, which converges. So $(v_{n})$ is bounded in $D^{1,p}(\mathbb{R}^N)$, and $v_{n}\to v$ almost everywhere. By the weak lower semicontinuity, $v\in D^{1,p}(\mathbb{R}^N)$ and
\[
S(N,p)^{p}\left(\int_{\mathbb{R}^N}|v|^{\frac{Np}{N-p}}dx\right)^{\frac{N-p}{N}}\leq\int_{\mathbb{R}^N}|\nabla v|^{p}dx
\leq\int_{\mathbb{H}^N}|\nabla_{\mathbb{H}}u^{\sharp}|^{p}dV_{\mathbb{H}}-\beta^{p}\int_{\mathbb{H}^N}(u^{\sharp})^{p}dV_{\mathbb{H}} .
\]
By our assumption, the right side is equal to $S(N,p)^{p}\left(\int_{\mathbb{H}^N}(u^{\sharp})^{\frac{Np}{N-p}}dV_{\mathbb{H}}\right)^{\frac{N-p}{N}}$. But \eqref{eq:norm} gives
\[
\int_{\mathbb{R}^N}|v|^{\frac{Np}{N-p}}dx>\int_{\mathbb{H}^N}(u^{\sharp})^{\frac{Np}{N-p}}dV_{\mathbb{H}}
\]
since $u^{\sharp}\not\equiv0$. This is a contradiction. The proof of Theorem \ref{thm:main} is complete. \qed

\medskip

\section{The range $p\geq3$ -- Proofs of Theorems \ref{thm:pge3} and \ref{thm:PS}}\label{sec:pge3}

In this section we prove Theorem \ref{thm:pge3} and Theorem \ref{thm:PS}. Here we keep the function and we change the radial variable. The variable $t=\sinh\rho$ is the radial variable of the hyperboloid model of $\mathbb{H}^N$. In this variable, the volume element for radial functions is $t^{N-1}(1+t^{2})^{-1/2}dt$, see for example \cite{Be15}. We use the argument of Liang, Lu, Wang and Yang \cite[Section 3]{LLWY}. For $p=N=3$, they got the spectral gap in this variable from the one-dimensional Hardy inequality. We follow their argument for all $N\geq3$ and $p\geq3$. For $p>3$ we also use the inequality $(1+t^{2})^{\frac{p-1}{2}}\geq1+t^{p-1}$.

For a radial function $u=u(\rho)$ on $\mathbb{H}^N$, we define a radial function on $\mathbb{R}^N$ by
\begin{equation}\label{eq:sub}
\tilde v(x):=u(\rho),\qquad |x|=t=\sinh\rho .
\end{equation}
Since $dt=\cosh\rho\,d\rho$ and $u'(\rho)=\tilde v'(t)\cosh\rho$, we have for all $p\geq1$ and $q\geq1$,
\begin{equation}\label{eq:subid}
\begin{gathered}
\int_{0}^{\infty}|u'|^{p}\sinh^{N-1}\rho\,d\rho=\int_{0}^{\infty}|\tilde v'|^{p}\,t^{N-1}(1+t^{2})^{\frac{p-1}{2}}dt,\\
\int_{0}^{\infty}|u|^{q}\sinh^{N-1}\rho\,d\rho=\int_{0}^{\infty}|\tilde v|^{q}\,\frac{t^{N-1}}{\sqrt{1+t^{2}}}\,dt .
\end{gathered}
\end{equation}
In particular, $\|\tilde v\|_{L^{q}(\mathbb{R}^N)}\geq\|u\|_{L^{q}(\mathbb{H}^N)}$ for all $q$, and the inequality is strict if $u\not\equiv0$. We will use the one-dimensional Hardy inequality \cite[Theorem 330]{HLP}: for $m>p\geq1$ and every Lipschitz function $w$ with compact support in $[0,\infty)$,
\begin{equation}\label{eq:hardy}
\int_{0}^{\infty}|w'|^{p}\,t^{m-1}dt\geq\left(\frac{m-p}{p}\right)^{p}\int_{0}^{\infty}|w|^{p}\,t^{m-p-1}dt .
\end{equation}
For the convenience of the reader we recall the proof. By integration by parts,
\[
\int_0^\infty|w|^{p}t^{m-p-1}dt=-\frac{p}{m-p}\int_0^\infty|w|^{p-2}w\,w'\,t^{m-p}dt,
\]
and the H\"older inequality gives \eqref{eq:hardy}.

\begin{proposition}\label{prop:sub}
Let $N\geq3$ and $p\geq3$. Let $u\geq0$ be a radial and Lipschitz function with compact support on $\mathbb{H}^N$, and let $\tilde v$ be as in \eqref{eq:sub}. Then
\begin{equation}\label{eq:keysub}
\int_{\mathbb{H}^N}|\nabla_{\mathbb{H}}u|^{p}dV_{\mathbb{H}}-\beta^{p}\int_{\mathbb{H}^N}u^{p}dV_{\mathbb{H}}
\geq\int_{\mathbb{R}^N}|\nabla\tilde v|^{p}dx+\beta^{p}\int_{\mathbb{H}^N}\frac{2}{e^{2\rho}-1}\,u^{p}dV_{\mathbb{H}} .
\end{equation}
Moreover, the last integral is equal to
\[
\beta^{p}|\mathbb{S}^{N-1}|\int_{0}^{\infty}\tilde v^{p}\,\frac{t^{N-2}}{\sqrt{1+t^{2}}\,\left(\sqrt{1+t^{2}}+t\right)}\,dt,
\]
and it is at least $\frac{\beta^{p}}{2}\int_{\mathbb{R}^N}\frac{\tilde v^{p}}{(1+|x|^{2})^{3/2}}\,dx$.
\end{proposition}

\begin{proof}
Let $a:=\frac{p-1}{2}\geq1$. Then
\[
(1+t^{2})^{a}\geq1+t^{2a}=1+t^{p-1},\qquad t\geq0,
\]
since $(x+y)^{a}\geq x^{a}+y^{a}$ for $x,y\geq0$ and $a\geq1$. So by \eqref{eq:subid},
\[
\int_{0}^{\infty}|u'|^{p}\sinh^{N-1}\rho\,d\rho
\geq\int_{0}^{\infty}|\tilde v'|^{p}\,t^{N-1}dt+\int_{0}^{\infty}|\tilde v'|^{p}\,t^{N+p-2}dt .
\]
We apply \eqref{eq:hardy} to the second integral with $m=N+p-1$. Then the constant is $\left(\frac{N-1}{p}\right)^{p}=\beta^{p}$, and the weight is $t^{N-2}$. So, by the second identity in \eqref{eq:subid},
\[
\int_{0}^{\infty}\left(|u'|^{p}-\beta^{p}u^{p}\right)\sinh^{N-1}\rho\,d\rho
\geq\int_{0}^{\infty}|\tilde v'|^{p}\,t^{N-1}dt
+\beta^{p}\int_{0}^{\infty}\tilde v^{p}\,t^{N-2}\left(1-\frac{t}{\sqrt{1+t^{2}}}\right)dt .
\]
Now
\[
1-\frac{t}{\sqrt{1+t^{2}}}=\frac{1}{\sqrt{1+t^{2}}\left(\sqrt{1+t^{2}}+t\right)}.
\]
In the variable $\rho$, this is $\frac{\cosh\rho-\sinh\rho}{\cosh\rho}$. So
\[
t^{N-2}\left(1-\frac{t}{\sqrt{1+t^{2}}}\right)dt=\sinh^{N-1}\rho\,\frac{e^{-\rho}}{\sinh\rho}\,d\rho
=\frac{2}{e^{2\rho}-1}\sinh^{N-1}\rho\,d\rho .
\]
We multiply by $|\mathbb{S}^{N-1}|$ and get \eqref{eq:keysub}. The last claim follows from $\sqrt{1+t^{2}}+t\leq2\sqrt{1+t^{2}}$ and $t^{N-2}\geq t^{N-1}(1+t^{2})^{-1/2}$.
\end{proof}

For $p>3$, we dropped the nonnegative term
\[
\int_0^\infty|\tilde v'|^{p}t^{N-1}\left[(1+t^{2})^{\frac{p-1}{2}}-1-t^{p-1}\right]dt
\]
in the proof. This term could give one more remainder by \eqref{eq:hardy}. For $p<3$, the inequality $(1+t^{2})^{a}\geq1+t^{p-1}$ fails for every $t>0$, so this proof does not work. So for $p<3$ we need Theorem \ref{thm:main}.

\begin{proof}[Proof of Theorem \ref{thm:pge3}]
The proof is the same as in Section \ref{sec:proof}. Let $u\geq0$ be a radial, nonincreasing and Lipschitz function with compact support. By Proposition \ref{prop:sub}, by the Sobolev inequality \eqref{eq:sobolev} for $\tilde v$, and by $\|\tilde v\|_{L^{\frac{Np}{N-p}}(\mathbb{R}^N)}\geq\|u\|_{L^{\frac{Np}{N-p}}(\mathbb{H}^N)}$, we have
\[
\int_{\mathbb{H}^N}|\nabla_{\mathbb{H}}u|^{p}dV_{\mathbb{H}}-\beta^{p}\int_{\mathbb{H}^N}u^{p}dV_{\mathbb{H}}
\geq S(N,p)^{p}\left(\int_{\mathbb{H}^N}|u|^{\frac{Np}{N-p}}dV_{\mathbb{H}}\right)^{\frac{N-p}{N}}+\beta^{p}\int_{\mathbb{H}^N}\frac{2}{e^{2\rho}-1}\,u^{p}dV_{\mathbb{H}} .
\]
This is \eqref{eq:mainW} for such $u$. The remainder is finite, since the weight is $O(\rho^{-1})$ near the origin. The weight $\frac{2}{e^{2\rho}-1}$ is positive and decreasing in $\rho$. So we can pass to all $u\in C_{c}^{\infty}(\mathbb{H}^N)$ by the P\'olya--Szeg\"o principle and the Hardy--Littlewood inequality, and then to all $u\in D^{1,p}(\mathbb{H}^N)$ by density and Fatou's lemma, as in Section \ref{sec:proof}. The sharpness of $\beta^{p}$ and of $S(N,p)^{p}$ is also proved as in Section \ref{sec:proof}.

Finally, assume that $u\not\equiv0$ gives equality in \eqref{eq:PS}. Then $u^{\sharp}$ also gives equality. We take radial, nonincreasing and Lipschitz functions $u_{n}$ with compact support such that $u_n\to u^{\sharp}$. The functions $\tilde v_{n}$ are bounded in $D^{1,p}(\mathbb{R}^N)$ and converge almost everywhere to $\tilde v$. So
\[
S(N,p)^{p}\|\tilde v\|^{p}_{L^{\frac{Np}{N-p}}(\mathbb{R}^N)}\leq\int_{\mathbb{R}^N}|\nabla\tilde v|^{p}dx\leq S(N,p)^{p}\|u^{\sharp}\|^{p}_{L^{\frac{Np}{N-p}}(\mathbb{H}^N)}.
\]
This is not possible, since $\|\tilde v\|_{L^{\frac{Np}{N-p}}(\mathbb{R}^N)}>\|u^{\sharp}\|_{L^{\frac{Np}{N-p}}(\mathbb{H}^N)}$.
\end{proof}

\begin{proof}[Proof of Theorem \ref{thm:PS}]
We drop the remainder terms in Theorem \ref{thm:main} and Theorem \ref{thm:pge3}. The sharpness of the two constants and the fact that $S(N,p)^p$ is not attained were proved together with these theorems. For $N=3$, the case (i) of Theorem \ref{thm:main} is the whole range. For $N=4$, Theorem \ref{thm:main} covers $\frac{1+\sqrt{13}}{2}\leq p\leq3$ and Theorem \ref{thm:pge3} covers $3\leq p<4$. For $N=5,6$, we have $\frac{1+\sqrt{4N-3}}{2}<3\leq\min\{5,\frac{N+1}{2}\}$. So Theorem \ref{thm:main} covers $\frac{1+\sqrt{4N-3}}{2}\leq p\leq 3$ and Theorem \ref{thm:pge3} covers $3\leq p<N$. For $N\geq7$, we have $\frac{1+\sqrt{4N-3}}{2}\geq3$, and Theorem \ref{thm:pge3} covers $3\leq p<N$.
\end{proof}

\medskip

\section{The full spectral gap for small $p$ -- Proof of Theorem \ref{thm:gap}}\label{sec:gap}

In this section we prove Theorem \ref{thm:gap}. Let $N\geq3$ and $1<p<\frac{1+\sqrt{4N-3}}{2}$. The upper bound for $p$ is the same as $p^{2}-p<N-1$, that is, $p-1<\beta$.

\subsection{A convexity estimate.}
The next lemma bounds from below the difference between the function $|t|^{p}$ at $t=a+b$ and its tangent line at $t=a$. Here $|a|^{p-2}a=0$ if $a=0$.

\begin{lemma}\label{lem:convex}
Let $a,b\in\mathbb{R}$.
\begin{itemize}
\item[(i)] If $1<p\leq2$ and $b\neq0$, then $|a+b|^{p}-|a|^{p}-p|a|^{p-2}a\,b\geq\frac{p(p-1)}{2}\,b^{2}\left(|a|+|b|\right)^{p-2}$.
\item[(ii)] If $p\geq2$, then $|a+b|^{p}-|a|^{p}-p|a|^{p-2}a\,b\geq p(p-1)4^{-p}|b|^{p}$.
\end{itemize}
\end{lemma}

\begin{proof}
We may assume $b\neq0$. Let $f(t)=|a+tb|^{p}$ for $0\leq t\leq1$. Since $p>1$, $f$ is $C^{1}$, and $f'$ is absolutely continuous with
\[
f''(t)=p(p-1)|a+tb|^{p-2}b^{2}
\]
for all $t$ with $a+tb\neq0$. There is at most one $t$ with $a+tb=0$, and $f''$ is integrable since $p-2>-1$. So
\[
|a+b|^{p}-|a|^{p}-p|a|^{p-2}a\,b=f(1)-f(0)-f'(0)=\int_{0}^{1}(1-t)f''(t)\,dt .
\]
(i) Since $|a+tb|\leq|a|+|b|$ and $p-2\leq0$, we have $f''(t)\geq p(p-1)b^{2}(|a|+|b|)^{p-2}$. Also $\int_0^1(1-t)\,dt=\frac12$.

(ii) The set of $t$ with $|a+tb|<\frac{|b|}{4}$ is the set of $t$ with $\left|t+\frac ab\right|<\frac14$. It is an interval of length $\frac12$. So the set of $t\in[0,\frac34]$ with $|a+tb|\geq\frac{|b|}{4}$ has measure at least $\frac14$. On this set we have $1-t\geq\frac14$ and $f''(t)\geq p(p-1)4^{2-p}|b|^{p}$, since $p-2\geq0$. Hence the integral is at least $\frac14\cdot\frac14\cdot p(p-1)4^{2-p}|b|^{p}$.
\end{proof}

\subsection{The ground state transform.}
Let
\[
\Lambda(r):=\left(\beta-p+1\right)\frac{\coth^{p-2}r}{\sinh^{2}r}+\beta\left(\coth^{p-2}r-1\right),\qquad r>0 .
\]

\begin{lemma}\label{lem:Lambda}
There is a constant $c>0$, which depends only on $N$ and $p$, such that
\[
\Lambda(r)\geq c\left(\coth^{p}r-1\right)>0\qquad\text{for all }r>0 .
\]
\end{lemma}

\begin{proof}
Let $x=\coth r>1$. Since $\frac{1}{\sinh^{2}r}=x^{2}-1$, we have $\Lambda(r)=h(x)$, where
\[
h(x)=(\beta-p+1)x^{p}+(p-1)x^{p-2}-\beta .
\]
Then $h(1)=0$ and
\[
h'(y)=p(\beta-p+1)y^{p-1}-(p-1)(2-p)y^{p-3}.
\]
Let $y\geq1$. If $p<2$, then $(p-1)(2-p)y^{p-3}\leq(p-1)(2-p)y^{p-1}$, and since $p\beta=N-1$,
\[
h'(y)\geq\left[p(\beta-p+1)-(p-1)(2-p)\right]y^{p-1}=(N+1-2p)y^{p-1},
\]
where $N+1-2p>0$ since $2p<4\leq N+1$. If $p\geq2$, the last term of $h'(y)$ is nonnegative, so $h'(y)\geq p(\beta-p+1)y^{p-1}$, where $\beta-p+1>0$ by our condition on $p$. In both cases $h'(y)\geq pc\,y^{p-1}$ for some $c>0$. We integrate from $1$ to $x$ and get $h(x)\geq c(x^{p}-1)$.
\end{proof}

\begin{proposition}\label{prop:ground}
Let $u\geq0$ be a radial and Lipschitz function with compact support on $\mathbb{H}^N$, and let $g(r)=u(r)\sinh^{\beta}r$. Then
\[
\int_{0}^{\infty}\left(|u'|^{p}-\beta^{p}u^{p}\right)\sinh^{N-1}r\,dr
=\beta^{p-1}\int_{0}^{\infty}\Lambda(r)\,g^{p}\,dr+\int_{0}^{\infty}\left(|a+b|^{p}-|a|^{p}-p|a|^{p-2}a\,b\right)dr ,
\]
where $a=-\beta\coth r\,g(r)$ and $b=g'(r)$.
\end{proposition}

\begin{proof}
We have $u'=\sinh^{-\beta}r\left(g'-\beta\coth r\,g\right)$ and $\beta p=N-1$. So
\[
|u'|^{p}\sinh^{N-1}r=|a+b|^{p},\qquad u^{p}\sinh^{N-1}r=g^{p}.
\]
Here $|a|^{p}=\beta^{p}\coth^{p}r\,g^{p}$ and
\[
p|a|^{p-2}a\,b=-p\beta^{p-1}\coth^{p-1}r\,g^{p-1}g'=-\beta^{p-1}\coth^{p-1}r\,(g^{p})' .
\]
We integrate the last term by parts. Since $(\coth r)'=-\frac{1}{\sinh^{2}r}$, we get
\[
\int_{0}^{\infty}p|a|^{p-2}a\,b\,dr=-(p-1)\beta^{p-1}\int_{0}^{\infty}\frac{\coth^{p-2}r}{\sinh^{2}r}\,g^{p}\,dr .
\]
The boundary terms are zero, since $u$ has compact support and $g(r)\leq Cr^{\beta}$ near $0$, so that $\coth^{p-1}r\,g^{p}\leq C r^{N-p}\to0$. All the integrals are finite. Indeed, $g\leq Cr^{\beta}$ and $|g'|\leq Cr^{\beta-1}$ near $0$, so all the integrands are bounded by $Cr^{N-1-p}$ near $0$, and $p<N$. Finally, as in the proof of Proposition \ref{prop:key},
\[
\begin{aligned}
\beta^{p}\left(\coth^{p}r-1\right)-(p-1)\beta^{p-1}\frac{\coth^{p-2}r}{\sinh^{2}r}
&=\beta^{p-1}\left[(\beta-p+1)\frac{\coth^{p-2}r}{\sinh^{2}r}+\beta\left(\coth^{p-2}r-1\right)\right]\\
&=\beta^{p-1}\Lambda(r).
\end{aligned}
\]
This proves the proposition.
\end{proof}

\subsection{Proof of Theorem \ref{thm:gap}}
By the P\'olya--Szeg\"o principle, as in Section \ref{sec:proof}, it is enough to prove \eqref{eq:gapC} for $u\geq0$ radial, nonincreasing and Lipschitz with compact support. Then \eqref{eq:gapC} holds for $u\in C_{c}^{\infty}(\mathbb{H}^N)$, and also for $u\in D^{1,p}(\mathbb{H}^N)$ by density, since the left side is continuous in $D^{1,p}(\mathbb{H}^N)$ by \eqref{eq:gap}, and the right side is continuous by the Sobolev inequality on $\mathbb{H}^N$.

So let $u$ be such a function, let $g=u\sinh^{\beta}r$, and let
\[
D:=\int_{0}^{\infty}\left(|u'|^{p}-\beta^{p}u^{p}\right)\sinh^{N-1}r\,dr .
\]
In the rest of the proof, $C$ denotes a positive constant which depends only on $N$ and $p$, and which can change from line to line. Both integrals on the right side of Proposition \ref{prop:ground} are nonnegative, by Lemma \ref{lem:Lambda} and Lemma \ref{lem:convex}. So
\begin{equation}\label{eq:gapD}
\int_{0}^{\infty}\Lambda(r)\,g^{p}\,dr\leq CD,\qquad
\int_{0}^{\infty}\left(|a+b|^{p}-|a|^{p}-p|a|^{p-2}a\,b\right)dr\leq D ,
\end{equation}
with $a=-\beta\coth r\,g$ and $b=g'$ as in Proposition \ref{prop:ground}. Also, since
\[
u^{\frac{Np}{N-p}}\sinh^{N-1}r=g^{\frac{Np}{N-p}}\sinh^{-\frac{(N-1)p}{N-p}}r ,
\]
it is enough to prove
\begin{equation}\label{eq:gapgoal}
\int_{0}^{\infty}g^{\frac{Np}{N-p}}\sinh^{-\frac{(N-1)p}{N-p}}r\,dr\leq C D^{\frac{N}{N-p}} .
\end{equation}
We will use two pointwise bounds for $|g'|$ in terms of
\[
e(r):=|a+b|^{p}-|a|^{p}-p|a|^{p-2}a\,b\geq0 .
\]

(F1) At the points where $|g'|\geq\beta\coth r\,g$, we have $|g'|^{p}\leq Ce(r)$. Indeed, for $p<2$ we have $|a|+|b|\leq2|g'|$ there, so Lemma \ref{lem:convex} (i) gives $e(r)\geq\frac{p(p-1)}{2}2^{p-2}|g'|^{p}$. For $p\geq2$, Lemma \ref{lem:convex} (ii) gives $|g'|^{p}\leq Ce(r)$ at every point.

(F2) Let $p<2$ and $r\geq1$. Then $|a|+|b|\leq\beta\coth1\,g+|g'|\leq C\max\{g,|g'|\}$. So by Lemma \ref{lem:convex} (i), at the points where $|g'|<g$ we have
\[
|g'|\leq C\,e(r)^{\frac12}g^{\frac{2-p}{2}},
\]
and at the points where $|g'|\geq g$ we have $|g'|\leq C\,e(r)^{\frac1p}$.

\emph{Near the origin.} On $(0,3]$ we have $\coth^{p}r-1\geq(1-\tanh^{p}3)\coth^{p}r$. So by Lemma \ref{lem:Lambda} and \eqref{eq:gapD},
\begin{equation}\label{eq:near}
\int_{0}^{3}\coth^{p}r\,g^{p}\,dr\leq CD.
\end{equation}
Let $\varphi(r)=1$ for $r\leq2$, $\varphi(r)=3-r$ for $2\leq r\leq3$, and $\varphi(r)=0$ for $r\geq3$. We apply the Sobolev inequality on $\mathbb{H}^N$ to $\varphi u$. We get
\[
\left(\int_{B_{2}}u^{\frac{Np}{N-p}}dV_{\mathbb{H}}\right)^{\frac{N-p}{N}}
\leq S(N,p)^{-p}\int_{\mathbb{H}^N}|\nabla_{\mathbb{H}}(\varphi u)|^{p}dV_{\mathbb{H}}
\leq C\int_{B_{3}}\left(|u'|^{p}+u^{p}\right)dV_{\mathbb{H}} .
\]
Now $\int_{B_{3}}u^{p}dV_{\mathbb{H}}=|\mathbb{S}^{N-1}|\int_{0}^{3}g^{p}dr$ and $g^{p}\leq\coth^{p}r\,g^{p}$. Also,
\[
\int_{B_{3}}|u'|^{p}dV_{\mathbb{H}}=|\mathbb{S}^{N-1}|\int_{0}^{3}|g'-\beta\coth r\,g|^{p}dr\leq C\int_{0}^{3}\left(|g'|^{p}+\coth^{p}r\,g^{p}\right)dr,
\]
and by (F1), $|g'|^{p}\leq Ce(r)+\beta^{p}\coth^{p}r\,g^{p}$ at every point. So by \eqref{eq:gapD} and \eqref{eq:near},
\[
\int_{0}^{2}g^{\frac{Np}{N-p}}\sinh^{-\frac{(N-1)p}{N-p}}r\,dr=\frac{1}{|\mathbb{S}^{N-1}|}\int_{B_{2}}u^{\frac{Np}{N-p}}dV_{\mathbb{H}}\leq C D^{\frac{N}{N-p}} .
\]

\emph{Far from the origin.} First we bound $g(2)$. Since $u$ is nonincreasing, for $1\leq t\leq2$ we have $g(t)\geq u(2)\sinh^{\beta}t\geq\left(\frac{\sinh1}{\sinh2}\right)^{\beta}g(2)$. Also $\Lambda(t)\geq c(\coth^{p}2-1)$ for $1\leq t\leq2$, by Lemma \ref{lem:Lambda}. So by \eqref{eq:gapD},
\[
g(2)^{p}\leq\left(\frac{\sinh2}{\sinh1}\right)^{N-1}\int_{1}^{2}g^{p}\,dt\leq CD.
\]
Now let $R>2$, and let $G$ be the maximum of $g$ on $[2,R]$. For $2\leq s\leq R$ we have $g(s)\leq g(2)+\int_{2}^{R}|g'|\,dr$.

If $p<2$, we use (F2), the Cauchy--Schwarz inequality, the H\"older inequality and \eqref{eq:gapD}. We get
\[
\int_{2}^{R}|g'|\,dr\leq C\,G^{\frac{2-p}{2}}(R-2)^{\frac12}D^{\frac12}+C(R-2)^{\frac{p-1}{p}}D^{\frac1p}.
\]
By the Young inequality $xy\leq\frac{2-p}{2}x^{\frac{2}{2-p}}+\frac p2y^{\frac2p}$ with $x=G^{\frac{2-p}{2}}$, the first term is at most $\frac{2-p}{2}G+C(R-2)^{\frac1p}D^{\frac1p}$. Since $(R-2)^{\frac1p}\leq R$ and $(R-2)^{\frac{p-1}{p}}\leq R$, we get
\[
G\leq g(2)+\frac{2-p}{2}G+CRD^{\frac1p}.
\]
Since $G$ is finite, this gives $G\leq\frac2p\,g(2)+CRD^{\frac1p}$. If $p\geq2$, then $|g'|^{p}\leq Ce(r)$ at every point by (F1), and the H\"older inequality and \eqref{eq:gapD} give directly $G\leq g(2)+C(R-2)^{\frac{p-1}{p}}D^{\frac1p}$. In both cases, by the bound for $g(2)$,
\[
g(R)\leq G\leq C\left(D^{\frac1p}+RD^{\frac1p}\right)\leq CRD^{\frac1p}\qquad\text{for all }R\geq2 .
\]
Hence
\[
\int_{2}^{\infty}g^{\frac{Np}{N-p}}\sinh^{-\frac{(N-1)p}{N-p}}r\,dr\leq C D^{\frac{N}{N-p}}\int_{2}^{\infty}r^{\frac{Np}{N-p}}\sinh^{-\frac{(N-1)p}{N-p}}r\,dr ,
\]
and the last integral is finite. Together with the bound near the origin, this proves \eqref{eq:gapgoal}, and so \eqref{eq:gapC}.

\emph{The best constant.} The test functions $u_{\varepsilon}$ in the proof of the sharpness of $S(N,p)^{p}$ in Section \ref{sec:proof} can be used for every $1<p<N$. They show that the best constant in \eqref{eq:gapC} is at most $S(N,p)^{p}$. When $p<\min\{2,\sqrt N\}$, Proposition \ref{prop:fail} gives a function $u$ for which \eqref{eq:gapC} with $C=S(N,p)^{p}$ is false. So in this case the best constant is strictly smaller than $S(N,p)^{p}$. This completes the proof of Theorem \ref{thm:gap}. \qed

\medskip

\section{The logarithmic Sobolev inequality -- Proofs of Theorems \ref{thm:logPS} and \ref{thm:log}}\label{sec:log}

In this section we prove Theorem \ref{thm:logPS} and Theorem \ref{thm:log}. In the cases of Theorem \ref{thm:logPS} we know that $R_{N}\geq0$. In the case (a) this follows from the paragraph after Lemma \ref{lem:mono}, in the case (b) from Lemma \ref{lem:mono}, and in the case (c) from Lemma \ref{lem:pos}. In the cases of Theorem \ref{thm:log} we also know that $R_{N}$ is nonincreasing, by Lemma \ref{lem:mono}. Note that the first three lines of \eqref{eq:chain} hold in all these cases, since they only use Proposition \ref{prop:key}.

For Corollary \ref{cor:log}, we replace \eqref{eq:chain} in the proof below by the comparison of Nguyen \cite[Theorem 2.2]{H18}, which gives
\[
\int_{\mathbb{H}^N}|\nabla_{\mathbb{H}}u^{\sharp}|^{p}dV_{\mathbb{H}}-\beta^{p}\int_{\mathbb{H}^N}(u^{\sharp})^{p}dV_{\mathbb{H}}\geq\int_{\mathbb{R}^N}|\nabla u^{\sharp}_{e}|^{p}dx
\]
for the Euclidean rearrangement $u^{\sharp}_{e}$ of $u$, see also \cite[Section 3]{DGLL}. The function $u^{\sharp}_{e}$ has the same $L^p$ norm and the same entropy as $u$. Everything else is the same.

\begin{proof}[Proof of Theorem \ref{thm:logPS} and Theorem \ref{thm:log}]
\emph{Radial functions.} Let $u\geq0$ be a radial, nonincreasing and Lipschitz function with compact support such that $\int_{\mathbb{H}^N}u^{p}dV_{\mathbb{H}}=1$. Let $v=um$ be as in \eqref{eq:map}. By \eqref{eq:transfer}, we have $\int_{\mathbb{R}^N}v^{p}dx=1$. Since $\ln v=\ln u+\beta\ln\frac{\sinh r}{r}$, \eqref{eq:transfer} also gives
\begin{equation}\label{eq:entmap}
\int_{\mathbb{R}^N}v^{p}\ln v\,dx=\int_{\mathbb{H}^N}u^{p}\ln u\,dV_{\mathbb{H}}+\beta\int_{\mathbb{H}^N}u^{p}\ln\frac{\sinh\rho}{\rho}\,dV_{\mathbb{H}} .
\end{equation}
The two integrals on the right are finite, since $u$ is bounded with compact support and $\ln\frac{\sinh\rho}{\rho}\leq\rho$. Now we apply \eqref{eq:logE} to $v$. Then we use the first three lines of \eqref{eq:chain} and the fact that the logarithm is increasing. We get
\begin{equation}\label{eq:logchain}
\begin{aligned}
\int_{\mathbb{H}^N}u^{p}\ln u\,dV_{\mathbb{H}}+\beta\int_{\mathbb{H}^N}u^{p}\ln\frac{\sinh\rho}{\rho}\,dV_{\mathbb{H}}
&\leq\frac{N}{p^{2}}\ln\left(\mathcal{L}_{N,p}\int_{\mathbb{R}^N}|\nabla v|^{p}dx\right)\\
&\leq\frac{N}{p^{2}}\ln\Bigg[\mathcal{L}_{N,p}\Bigg(\int_{\mathbb{H}^N}|\nabla_{\mathbb{H}}u|^{p}dV_{\mathbb{H}}-\beta^{p}\\
&\qquad\qquad-\beta^{p-1}\int_{\mathbb{H}^N}R_{N}(\rho)\,u^{p}dV_{\mathbb{H}}\Bigg)\Bigg].
\end{aligned}
\end{equation}
For $N=3$ and $p=3$, we use that \eqref{eq:logE} holds for every $p>1$ \cite{IG}. This is the statement of Theorem \ref{thm:log} for radial nonincreasing functions. If we drop the second term on the left, which is nonnegative, and the term with $R_{N}$, which is also nonnegative, we get Theorem \ref{thm:logPS} for such functions.

\emph{All functions.} Let $u\in C_{c}^{\infty}(\mathbb{H}^N)$ with $\int_{\mathbb{H}^N}|u|^{p}dV_{\mathbb{H}}=1$, and let $u^{\sharp}$ be as in Section \ref{sec:proof}. By the layer cake formula, applied to the positive part and to the negative part of the function $s\mapsto s^p\ln s$, we have
\[
\int_{\mathbb{H}^N}|u|^{p}\ln|u|\,dV_{\mathbb{H}}=\int_{\mathbb{H}^N}(u^{\sharp})^{p}\ln u^{\sharp}\,dV_{\mathbb{H}},
\]
see also \cite[Section 3]{DGLL}. For Theorem \ref{thm:logPS}, we apply the inequality for radial functions to $u^{\sharp}$, and we use the P\'olya--Szeg\"o principle and the fact that the logarithm is increasing. For Theorem \ref{thm:log}, we apply \eqref{eq:logchain} to $u^{\sharp}$, we drop the term with $\ln\frac{\sinh\rho}{\rho}$, and we use the P\'olya--Szeg\"o principle and the Hardy--Littlewood inequality as in Section \ref{sec:proof}. This gives \eqref{eq:log} for $u$.

\emph{Functions in $W^{1,p}(\mathbb{H}^N)$.} Let $u\in W^{1,p}(\mathbb{H}^N)$ with $\int_{\mathbb{H}^N}|u|^{p}dV_{\mathbb{H}}=1$. The positive part of $\int_{\mathbb{H}^N}|u|^p\ln|u|\,dV_{\mathbb{H}}$ is finite. Indeed, for small $\delta>0$, $s^{p}\ln s\leq C_{\delta}s^{p+\delta}$ for $s\geq1$, and $u\in L^{p+\delta}(\mathbb{H}^N)$ by \eqref{eq:gap} and the Sobolev inequality on $\mathbb{H}^N$ (for $p=N$, by the Moser--Trudinger inequality). If the negative part is infinite, then the left side of \eqref{eq:log} is $-\infty$ and there is nothing to prove. In the other case, fix a nonincreasing function $\varphi\in C^\infty([0,\infty))$ with $\varphi=1$ on $[0,1]$ and $\varphi=0$ on $[2,\infty)$, and let $\chi_n(x)=\varphi(\rho(x)/n)$. Then $0\leq\chi_n\leq1$, $\chi_n$ is nondecreasing in $n$, $\chi_n\to1$, and $|\nabla_{\mathbb{H}}\chi_n|\leq\frac{C}{n}$. Let $u_{n}=c_{n}\chi_{n}u$, where $c_{n}\to1$ is chosen so that $\int_{\mathbb{H}^N}|u_n|^{p}dV_{\mathbb{H}}=1$. Then we mollify each $u_{n}$ and normalize again in $L^p$. With this approximation, all terms in \eqref{eq:log} converge. For the entropy term, this follows from the dominated convergence theorem, since $s^{p}|\ln s|$ is increasing on $(0,e^{-1/p})$ and continuous on $(0,\infty)$. For the remainder term, it follows from the monotone convergence of $\chi_{n}|u|$ and from $R_{N}\in L^{1}_{\mathrm{loc}}(\mathbb{H}^N)$. For the Dirichlet integral, it follows from $\nabla_{\mathbb{H}}(\chi_nu)=\chi_n\nabla_{\mathbb{H}}u+u\nabla_{\mathbb{H}}\chi_n$, $|\nabla_{\mathbb{H}}\chi_n|\leq\frac Cn$ and $u\in W^{1,p}(\mathbb{H}^N)$. Finally, the right side of \eqref{eq:log} is continuous under this convergence, since its argument stays away from $0$ by \eqref{eq:chain}.

\emph{The constant $\beta^{p}$ is sharp.} If \eqref{eq:logPS} holds with some $c>\beta^{p}$ in place of $\beta^{p}$, then
\[
\int_{\mathbb{H}^N}|\nabla_{\mathbb{H}}u|^{p}dV_{\mathbb{H}}-c>0
\]
for all $u\in C_{c}^{\infty}(\mathbb{H}^N)$ with $\int_{\mathbb{H}^N}|u|^{p}dV_{\mathbb{H}}=1$. This is not possible, since the constant in \eqref{eq:gap} is sharp.

\emph{The constant $\mathcal{L}_{N,p}$ is sharp.} Let $w\in C_{c}^{\infty}(\mathbb{R}^N)$ with $\int_{\mathbb{R}^N}|w|^{p}dy=1$, and let $u_{\varepsilon}(\exp_{0}(\varepsilon y))=\varepsilon^{-N/p}w(y)$ as in Section \ref{sec:proof}. Then $\int_{\mathbb{H}^N}|u_{\varepsilon}|^{p}dV_{\mathbb{H}}=1+O(\varepsilon^{2})$. After we normalize $u_\varepsilon$ in $L^p$, we have
\[
\int_{\mathbb{H}^N}|u_\varepsilon|^p\ln|u_\varepsilon|\,dV_{\mathbb{H}}=\int_{\mathbb{R}^N}|w|^p\ln|w|\,dy-\frac Np\ln\varepsilon+o(1)
\]
and
\[
\int_{\mathbb{H}^N}|\nabla_{\mathbb{H}}u_{\varepsilon}|^{p}dV_{\mathbb{H}}-\beta^{p}=\varepsilon^{-p}\left(\int_{\mathbb{R}^N}|\nabla w|^{p}dy+o(1)\right).
\]
So if \eqref{eq:logPS} holds with some $L<\mathcal{L}_{N,p}$ in place of $\mathcal{L}_{N,p}$, then as $\varepsilon\to0$ we get
\[
\int_{\mathbb{R}^N}|w|^p\ln|w|\,dy\leq\frac{N}{p^{2}}\ln\left(L\int_{\mathbb{R}^N}|\nabla w|^{p}dy\right)
\]
for all such $w$. This is not possible, since $\mathcal{L}_{N,p}$ is sharp in \eqref{eq:logE}.

\emph{The constant $\mathcal{L}_{N,p}$ is not attained.} Assume that $u\in W^{1,p}(\mathbb{H}^N)$ with $\int_{\mathbb{H}^N}|u|^{p}dV_{\mathbb{H}}=1$ gives equality in \eqref{eq:logPS}. Then $u^{\sharp}$ also gives equality, since the left side does not change and the right side does not increase. Let $v=u^{\sharp}m$. As in Section \ref{sec:proof}, we have $v\in W^{1,p}(\mathbb{R}^N)$, $\int_{\mathbb{R}^N}v^{p}dx=1$, and
\[
\int_{\mathbb{R}^N}|\nabla v|^{p}dx\leq\int_{\mathbb{H}^N}|\nabla_{\mathbb{H}}u^{\sharp}|^{p}dV_{\mathbb{H}}-\beta^{p}
\]
by the weak lower semicontinuity. By \eqref{eq:logE}, by this bound and by the equality for $u^{\sharp}$, we get
\[
\int_{\mathbb{R}^N}v^{p}\ln v\,dx\leq\int_{\mathbb{H}^N}(u^{\sharp})^{p}\ln u^{\sharp}\,dV_{\mathbb{H}}.
\]
But \eqref{eq:entmap} is an identity between the integrands, by \eqref{eq:transfer}. So it holds for $u^{\sharp}$, with both sides in $(-\infty,\infty]$, and its last term is positive since $u^{\sharp}\not\equiv0$. This is a contradiction.
\end{proof}

\medskip

\section{Hardy--Trudinger--Moser inequalities -- Proofs of Theorems \ref{thm:HTMN} and \ref{thm:HTM3}}\label{sec:TM}

In this section we study the endpoint $p=N\geq3$ and we prove Theorem \ref{thm:HTMN} and Theorem \ref{thm:HTM3}. When $p=N$, Proposition \ref{prop:sub} can be applied. Since $\beta=\frac{N-1}{N}$ in this case, it gives
\begin{equation}\label{eq:TMkey}
\begin{aligned}
\int_{\mathbb{H}^N}|\nabla_{\mathbb{H}}u|^{N}dV_{\mathbb{H}}-\left(\frac{N-1}{N}\right)^{N}\int_{\mathbb{H}^N}|u|^{N}dV_{\mathbb{H}}
\geq\ &\int_{\mathbb{R}^N}|\nabla\tilde v|^{N}dx\\
&+\frac12\left(\frac{N-1}{N}\right)^{N}\int_{\mathbb{R}^N}\frac{|\tilde v|^{N}}{(1+|x|^{2})^{3/2}}\,dx
\end{aligned}
\end{equation}
for all nonnegative, radial and Lipschitz functions $u$ with compact support, where $\tilde v(x)=u(\arsinh|x|)$. This is a version of \cite[Lemma 2.2]{LLWY}, which is for $N\geq4$, and of \cite[(3.4)]{LLWY}, which is for $N=3$, with the explicit weight $(1+|x|^{2})^{-3/2}$ and with a proof which is the same for all $N\geq3$. For $N=3$, the constant $\frac12\left(\frac23\right)^3=\frac{4}{27}$ is the same as in \cite[(3.4)]{LLWY}. In \cite{LLWY}, the Euclidean variable is the one of Nguyen for $N\geq4$, and it is $t=\sinh\rho$ for $N=3$. For $N\geq4$, the lemma in \cite{LLWY} is proved by the asymptotics of the function
\[
F(\rho)=\sinh^{N(N-1)}\rho-s^{N(N-1)}-\left(\frac{N-1}{N}\right)^{N}s^{N^{2}},
\]
see \cite[Lemma 2.1]{LLWY}, with different weights for $N=4$ and $N\geq5$. It is remarked in \cite{LLWY} that this way does not work for $N=3$. Proposition \ref{prop:sub} shows that the way of the three-dimensional case works for all $N\geq3$.

The rest of the proof of \eqref{eq:HTM} in \cite[Section 2]{LLWY} only uses the form of \eqref{eq:TMkey}, that is, a Dirichlet integral on $\mathbb{R}^N$ plus a weighted $L^{N}$ norm, where the weight is bounded from below by a multiple of $(1+|x|^{2})^{-\alpha}$ for some $\alpha>0$. It also uses the Hardy--Trudinger--Moser inequality on $\mathbb{R}^N$ \cite[Theorem 1.8]{LLWY}:
\begin{equation}\label{eq:LLWY18}
\int_{\mathbb{R}^N}|\nabla w|^{N}dx+\tau\int_{\mathbb{R}^N}\frac{|w|^{N}}{(1+|x|^{2})^{\alpha}}dx\leq1
\quad\Longrightarrow\quad
\int_{\mathbb{R}^N}\frac{\Phi_{N-1}\left(\beta_{N}|w|^{\frac{N}{N-1}}\right)}{(1+|x|^{2})^{\alpha}}dx\leq C(N,\alpha,\tau),
\end{equation}
and the Hardy--Sobolev--Maz'ya inequality on $\mathbb{H}^N$ \cite{Maz85, MST13}. The last inequality says that for every $q>N$ there is $S_{q}>0$ such that
\begin{equation}\label{eq:HSM}
\int_{\mathbb{H}^N}|\nabla_{\mathbb{H}}u|^{N}dV_{\mathbb{H}}-\left(\frac{N-1}{N}\right)^{N}\int_{\mathbb{H}^N}|u|^{N}dV_{\mathbb{H}}
\geq S_{q}\left(\int_{\mathbb{H}^N}|u|^{q}dV_{\mathbb{H}}\right)^{N/q},\qquad u\in W^{1,N}_{0}(\mathbb{H}^N).
\end{equation}
We follow this argument, and we keep the remainder term in \eqref{eq:keysub}.

\begin{proof}[Proof of Theorem \ref{thm:HTMN}]
Let
\[
A(u)=\int_{\mathbb{H}^N}|\nabla_{\mathbb{H}}u|^{N}dV_{\mathbb{H}}-\beta^{N}\int_{\mathbb{H}^N}|u|^{N}dV_{\mathbb{H}},\qquad
\mathcal{W}(u)=\beta^{N}\int_{\mathbb{H}^N}\frac{2}{e^{2\rho}-1}|u|^{N}dV_{\mathbb{H}},
\]
where $\beta=\frac{N-1}{N}$. Then the constraint \eqref{eq:constraintN} is $A(u)-\lambda \mathcal{W}(u)\leq1$.

\emph{Symmetrization.} Let $u^{\sharp}$ be the symmetric decreasing rearrangement of $|u|$. Then $A(u^{\sharp})\leq A(u)$ by the P\'olya--Szeg\"o principle. Also, $\mathcal{W}(u^{\sharp})\geq \mathcal{W}(u)$ by the Hardy--Littlewood inequality, since the weight $\frac{2}{e^{2\rho}-1}$ is positive and decreasing. The integrals in the conclusion do not change. So $u^{\sharp}$ also satisfies \eqref{eq:constraintN}, and it is enough to prove the theorem for $u\geq0$ radial and nonincreasing. By a standard approximation, we can also assume that $u$ is Lipschitz with compact support.

\emph{The Euclidean function.} Let $\tilde v(x)=u(\arsinh|x|)$ be as in \eqref{eq:sub}. By Proposition \ref{prop:sub} with $p=N$,
\begin{equation}\label{eq:AN}
A(u)\geq\int_{\mathbb{R}^N}|\nabla \tilde v|^{N}dx+\mathcal{W}(u),\qquad
\mathcal{W}(u)\geq\frac{\beta^{N}}{2}\int_{\mathbb{R}^N}\frac{\tilde v^{N}}{(1+|x|^{2})^{3/2}}dx .
\end{equation}
In particular, $\mathcal{W}(u)\leq A(u)$. So $(1-\lambda)A(u)\leq A(u)-\lambda \mathcal{W}(u)\leq1$, that is,
\begin{equation}\label{eq:ANbound}
A(u)\leq\frac{1}{1-\lambda}.
\end{equation}
By \eqref{eq:AN} again,
\[
\int_{\mathbb{R}^N}|\nabla \tilde v|^{N}dx+\tau\int_{\mathbb{R}^N}\frac{\tilde v^{N}}{(1+|x|^{2})^{3/2}}dx\leq \int_{\mathbb{R}^N}|\nabla \tilde v|^{N}dx+(1-\lambda)\mathcal{W}(u)\leq A(u)-\lambda \mathcal{W}(u)\leq 1,
\]
where $\tau:=\frac{1-\lambda}{2}\beta^{N}>0$. So \eqref{eq:LLWY18} with $\alpha=\frac32$ gives
\begin{equation}\label{eq:PhiN1}
\int_{\mathbb{R}^N}\frac{\Phi_{N-1}\left(\beta_{N}\tilde v^{\frac{N}{N-1}}\right)}{(1+|x|^{2})^{3/2}}dx\leq C_{0}(N,\lambda).
\end{equation}

\emph{The set where $u>1$.} Since $u$ is radial and nonincreasing, the set $\{u>1\}$ is a ball $B_{\rho_{0}}$ with center at the origin. Let $q=\frac{N^{2}}{N-1}>N$. By \eqref{eq:HSM} and by \eqref{eq:ANbound},
\[
|B_{\rho_{0}}|\leq\int_{B_{\rho_{0}}}u^{q}dV_{\mathbb{H}}\leq\left(\frac{A(u)}{S_{q}}\right)^{q/N}
\leq\left((1-\lambda)S_{q}\right)^{-q/N}.
\]
So $\rho_{0}\leq\rho_{1}$ for some $\rho_{1}=\rho_1(N,\lambda)<\infty$. Let $t_1=\sinh\rho_1$. Since $\Phi_{N}\leq\Phi_{N-1}$ and $u(\rho)=\tilde v(t)$ with $t=\sinh\rho$, we have by \eqref{eq:subid} and \eqref{eq:PhiN1},
\[
\begin{aligned}
\int_{B_{\rho_{0}}}\Phi_{N}\left(\beta_{N}u^{\frac{N}{N-1}}\right)dV_{\mathbb{H}}
&\leq|\mathbb{S}^{N-1}|\int_{0}^{t_{1}}\Phi_{N-1}\left(\beta_{N}\tilde v^{\frac{N}{N-1}}\right)\frac{t^{N-1}}{\sqrt{1+t^{2}}}\,dt\\
&\leq(1+t_{1}^{2})\,|\mathbb{S}^{N-1}|\int_{0}^{\infty}\Phi_{N-1}\left(\beta_{N}\tilde v^{\frac{N}{N-1}}\right)\frac{t^{N-1}}{(1+t^{2})^{3/2}}\,dt\\
&=(1+t_{1}^{2})\int_{\mathbb{R}^N}\frac{\Phi_{N-1}\left(\beta_{N}\tilde v^{\frac{N}{N-1}}\right)}{(1+|x|^{2})^{3/2}}dx\leq(1+t_{1}^{2})\,C_{0}.
\end{aligned}
\]

\emph{The set where $u\leq1$.} On this set we have
\[
\Phi_{N}\left(\beta_{N}u^{\frac{N}{N-1}}\right)=\sum_{n\geq N}\frac{\beta_{N}^{n}u^{\frac{nN}{N-1}}}{n!}\leq e^{\beta_{N}}u^{q}.
\]
By \eqref{eq:HSM} and \eqref{eq:ANbound}, $\int_{\mathbb{H}^N}u^{q}dV_{\mathbb{H}}\leq\left((1-\lambda)S_{q}\right)^{-q/N}$. The sum of the two bounds proves the first inequality.

\emph{The second inequality.} Since $dV_{\mathbb{H}}\geq dx$ on $B^N$, we have
\[
\int_{B^N}e^{\beta_{N}|u|^{\frac{N}{N-1}}}dx\leq\int_{\mathbb{H}^N}\Phi_{N}\left(\beta_{N}|u|^{\frac{N}{N-1}}\right)dV_{\mathbb{H}}+\sum_{k=0}^{N-1}\frac{\beta_N^k}{k!}\int_{B^N}|u|^{\frac{kN}{N-1}}dx .
\]
For $k=0$ the last integral is the volume of $B^N$. For $1\leq k\leq N-1$, we take $s>N-1$, so that $\frac{kNs}{N-1}>N$. By the H\"older inequality, $dx\leq dV_{\mathbb{H}}$, \eqref{eq:HSM} and \eqref{eq:ANbound},
\[
\int_{B^N}|u|^{\frac{kN}{N-1}}dx\leq\left(\int_{\mathbb{H}^N}|u|^{\frac{kNs}{N-1}}dV_{\mathbb{H}}\right)^{1/s}|B^N|^{1-\frac1s}\leq C(N,\lambda),
\]
as in \cite[Section 3]{LLWY}. This proves the second inequality.
\end{proof}

The proof of Theorem \ref{thm:HTM3} follows the same steps, with the weight $R_{3}$ in place of $\frac{2}{e^{2\rho}-1}$ and with the function $v$ of \eqref{eq:map} in place of $\tilde v$. We give the details, since the weight $R_3$ needs a different lower bound and a different comparison of $u$ and $v$.

\begin{proof}[Proof of Theorem \ref{thm:HTM3}]
Let
\[
A(u)=\int_{\mathbb{H}^3}|\nabla_{\mathbb{H}}u|^{3}dV_{\mathbb{H}}-\frac{8}{27}\int_{\mathbb{H}^3}|u|^{3}dV_{\mathbb{H}},\qquad
\mathcal{R}(u)=\frac49\int_{\mathbb{H}^3}R_{3}(\rho)|u|^{3}dV_{\mathbb{H}}.
\]
Then the constraint \eqref{eq:constraint3} is $A(u)-\lambda \mathcal{R}(u)\leq1$.

\emph{Symmetrization.} Let $u^{\sharp}$ be the symmetric decreasing rearrangement of $|u|$. Then $A(u^{\sharp})\leq A(u)$ by the P\'olya--Szeg\"o principle. Also, $\mathcal{R}(u^{\sharp})\geq \mathcal{R}(u)$ by the Hardy--Littlewood inequality, since $R_{3}$ is positive and decreasing by Lemma \ref{lem:mono} with $p=3$. The integral $\int_{\mathbb{H}^3}\Phi_{3}(\beta_{3}|u|^{3/2})dV_{\mathbb{H}}$ does not change. So $u^{\sharp}$ also satisfies \eqref{eq:constraint3}, and it is enough to prove the theorem for $u\geq0$ radial and nonincreasing. By a standard approximation, we can also assume that $u$ is Lipschitz with compact support.

\emph{The Euclidean function.} Let $v=u\left(\frac{\sinh r}{r}\right)^{2/3}$ be as in \eqref{eq:map}. By Proposition \ref{prop:key} with $p=3$ and by \eqref{eq:transfer},
\begin{equation}\label{eq:A3}
A(u)\geq\int_{\mathbb{R}^3}|\nabla v|^{3}dx+\mathcal{R}(u),\qquad
\mathcal{R}(u)=\frac49\,|\mathbb{S}^{2}|\int_{0}^{\infty}R_{3}(r)\,v^{3}\,r^{2}dr .
\end{equation}
In particular, $\mathcal{R}(u)\leq A(u)$. So $(1-\lambda)A(u)\leq A(u)-\lambda \mathcal{R}(u)\leq1$, that is,
\begin{equation}\label{eq:Abound}
A(u)\leq\frac{1}{1-\lambda}.
\end{equation}
By \eqref{eq:A3} again,
\[
\int_{\mathbb{R}^3}|\nabla v|^{3}dx+(1-\lambda)\mathcal{R}(u)\leq A(u)-\lambda \mathcal{R}(u)\leq1.
\]
The function $R_{3}$ is continuous and positive on $(0,\infty)$, and by Lemma \ref{lem:mono},
\[
\lim_{r\to0^{+}}r\,R_{3}(r)=\frac23,\qquad \lim_{r\to\infty}r^{3}R_{3}(r)=\frac43.
\]
So the number $c_{3}:=\inf_{r>0}(1+r^{2})^{3/2}R_{3}(r)$ is positive. Therefore
\[
\int_{\mathbb{R}^3}|\nabla v|^{3}dx+\tau\int_{\mathbb{R}^3}\frac{v^{3}}{(1+|x|^{2})^{3/2}}dx\leq1,
\qquad \tau:=(1-\lambda)\frac49c_{3},
\]
and \eqref{eq:LLWY18} with $N=3$ and $\alpha=\frac32$ gives
\begin{equation}\label{eq:Phi2}
\int_{\mathbb{R}^3}\frac{\Phi_{2}\left(\beta_{3}v^{3/2}\right)}{(1+|x|^{2})^{3/2}}dx\leq C_{0}(\lambda).
\end{equation}

\emph{The set where $u>1$.} Since $u$ is radial and nonincreasing, the set $\{u>1\}$ is a ball $B_{\rho_{0}}$ with center at the origin. By \eqref{eq:HSM} with $q=\frac92$ and by \eqref{eq:Abound},
\[
|B_{\rho_{0}}|\leq\int_{B_{\rho_{0}}}u^{9/2}dV_{\mathbb{H}}\leq\left(\frac{A(u)}{S_{9/2}}\right)^{3/2}
\leq\left((1-\lambda)S_{9/2}\right)^{-3/2}.
\]
Since $|B_{\rho}|=|\mathbb{S}^{2}|\int_{0}^{\rho}\sinh^{2}t\,dt$ goes to infinity as $\rho\to\infty$, we get $\rho_{0}\leq\rho_{1}(\lambda)$ for some $\rho_{1}(\lambda)<\infty$. On $B_{\rho_{0}}$ we have $u\leq v$, since $r\leq\sinh r$. Also, $\sinh^{2}r\leq C_{1}\,r^{2}(1+r^{2})^{-3/2}$ for $r\leq\rho_{1}$, with $C_{1}=C_{1}(\lambda)$, and $\Phi_{3}\leq\Phi_{2}$. So by \eqref{eq:Phi2},
\[
\int_{B_{\rho_{0}}}\Phi_{3}\left(\beta_{3}u^{3/2}\right)dV_{\mathbb{H}}
\leq|\mathbb{S}^{2}|\int_{0}^{\rho_{1}}\Phi_{2}\left(\beta_{3}v^{3/2}\right)\sinh^{2}r\,dr
\leq C_{1}\int_{\mathbb{R}^3}\frac{\Phi_{2}\left(\beta_{3}v^{3/2}\right)}{(1+|x|^{2})^{3/2}}dx\leq C_{1}C_{0}.
\]

\emph{The set where $u\leq1$.} On this set we have
\[
\Phi_{3}\left(\beta_{3}u^{3/2}\right)=\sum_{n\geq3}\frac{\beta_{3}^{n}u^{3n/2}}{n!}\leq e^{\beta_{3}}u^{9/2}.
\]
By \eqref{eq:HSM} and \eqref{eq:Abound}, $\int_{\mathbb{H}^3}u^{9/2}dV_{\mathbb{H}}\leq\left((1-\lambda)S_{9/2}\right)^{-3/2}$. The sum of the two bounds proves the first inequality. The second inequality follows from the first one, \eqref{eq:HSM} and \eqref{eq:Abound}, exactly as in the proof of Theorem \ref{thm:HTMN}.
\end{proof}

\medskip

\section{The case $p<2$}\label{sec:concl}

In this section we show that \eqref{eq:PS} and \eqref{eq:logPS} fail for some $p<2$. At the end we recall the weaker inequalities of Beckner and Nguyen, which hold for all $1\leq p<N$.

\medskip
Theorem \ref{thm:PS} starts at $p=2$, and this is not only a limit of our method. The next proposition shows that \eqref{eq:PS} fails below $p=2$ in a large range, in every dimension. For $p=1$, we read \eqref{eq:PS} for the characteristic function of a set $E$, with the perimeter $P(E)$ on the left side and the isoperimetric constant $S(N,1)=N|B|^{1/N}$ on the right side, where $|B|$ is the volume of the unit ball in $\mathbb{R}^N$.

\begin{proposition}\label{prop:fail}
Let $N\geq3$. The inequality \eqref{eq:PS} fails for $p=1$ and for every $1<p<\min\{2,\sqrt N\}$.
\end{proposition}

\begin{proof}
\emph{The case $p=1$.} For the geodesic ball $B_{R}$, we have
\[
P(B_{R})=|\mathbb{S}^{N-1}|\sinh^{N-1}R,\qquad
|B_{R}|=|\mathbb{S}^{N-1}|\int_{0}^{R}\sinh^{N-1}t\,dt.
\]
So
\[
\frac{P(B_{R})-(N-1)|B_{R}|}{(N-1)|\mathbb{S}^{N-1}|}=\int_{0}^{R}\sinh^{N-2}t\,(\cosh t-\sinh t)\,dt
=\int_{0}^{R}\sinh^{N-2}t\,e^{-t}dt .
\]
This grows like $R$ for $N=3$ and like $e^{(N-3)R}$ for $N\geq4$. The right side of \eqref{eq:PS}, which is $S(N,1)|B_{R}|^{\frac{N-1}{N}}$, grows like $e^{\frac{(N-1)^{2}}{N}R}$. Since $\frac{(N-1)^{2}}{N}>N-3$ for every $N$, \eqref{eq:PS} fails for large balls.

\emph{The case $1<p<\min\{2,\sqrt N\}$.} Let
\[
U(y)=\left(1+|y|^{\frac{p}{p-1}}\right)^{-\frac{N-p}{p}}
\]
be the Aubin--Talenti function, so that
\[
\int_{\mathbb{R}^N}|\nabla U|^{p}dy=S(N,p)^{p}\left(\int_{\mathbb{R}^N}|U|^{\frac{Np}{N-p}}dy\right)^{\frac{N-p}{N}}.
\]
At infinity, $|U|^{p}$ behaves like $|y|^{-\frac{p(N-p)}{p-1}}$ and $|\nabla U|^{p}$ behaves like $|y|^{-\frac{p(N-1)}{p-1}}$. Since $p^{2}<N$, the function $|U|^p$ is integrable on $\mathbb{R}^N$, that is, $U\in L^{p}(\mathbb{R}^N)$. Fix a function $\varphi\in C_{c}^{\infty}([0,1))$ with $\varphi=1$ on $[0,\frac12]$, and define on $\mathbb{H}^N$, in geodesic polar coordinates,
\[
u_{\varepsilon}(\rho)=\varepsilon^{-\frac{N-p}{p}}\,U(\rho/\varepsilon)\,\varphi(\rho),\qquad0<\varepsilon<1 .
\]
Let $\gamma=\frac{N-p}{p-1}$. In the variable $y=\rho/\varepsilon$, the volume element is $\left(\frac{\sinh(\varepsilon|y|)}{\varepsilon|y|}\right)^{N-1}dy$, and $\left(\frac{\sinh x}{x}\right)^{N-1}=1+O(x^{2})$ for $x\leq1$. The region $\rho\geq\frac12$, that is, $|y|\geq\frac{1}{2\varepsilon}$, gives $O(\varepsilon^{\gamma})$ in each of the three integrals
\[
\int_{\mathbb{H}^N}|\nabla_{\mathbb{H}}u_{\varepsilon}|^{p}dV_{\mathbb{H}},\qquad
\int_{\mathbb{H}^N}|u_{\varepsilon}|^{p}dV_{\mathbb{H}},\qquad
\int_{\mathbb{H}^N}|u_{\varepsilon}|^{\frac{Np}{N-p}}dV_{\mathbb{H}},
\]
by the decay of $U$ and $\nabla U$. For the $L^{p}$ integral, this part is $\varepsilon^{p}\,O\left(\varepsilon^{\frac{N-p^{2}}{p-1}}\right)=O(\varepsilon^{\gamma})$, which is $o(\varepsilon^{p})$ since $p<\gamma$. The term $O(\varepsilon^{2}|y|^{2})$ from the curvature on the region $|y|\leq\frac{1}{2\varepsilon}$ gives $O(\varepsilon^{2})$ when the second moment of the integrand is finite, and $O(\varepsilon^{\gamma}|\ln\varepsilon|)$ in the other case. Hence
\[
\begin{gathered}
\int_{\mathbb{H}^N}|\nabla_{\mathbb{H}}u_{\varepsilon}|^{p}dV_{\mathbb{H}}-S(N,p)^{p}\left(\int_{\mathbb{H}^N}|u_{\varepsilon}|^{\frac{Np}{N-p}}dV_{\mathbb{H}}\right)^{\frac{N-p}{N}}
=O\left(\varepsilon^{\min\{2,\gamma\}}|\ln\varepsilon|\right),\\
\beta^{p}\int_{\mathbb{H}^N}|u_{\varepsilon}|^{p}dV_{\mathbb{H}}=\beta^{p}\varepsilon^{p}\left(\int_{\mathbb{R}^N}|U|^{p}dy\right)\left(1+o(1)\right).
\end{gathered}
\]
In the first line we used the equality case of the Sobolev inequality for $U$. Since $p<2$ and $p<\gamma$ (the second condition is $p^{2}<N$), the left side of \eqref{eq:PS} minus the right side is
\[
-\beta^{p}\varepsilon^{p}\left(\int_{\mathbb{R}^N}|U|^{p}dy\right)\left(1+o(1)\right)<0
\]
for small $\varepsilon$.
\end{proof}

The logarithmic Sobolev inequality fails for every $1<p<2$ in every dimension, without the condition $p^{2}<N$. The reason is that the extremal functions in $\mathbb{R}^N$ decay fast.

\begin{proposition}\label{prop:faillog}
Let $N\geq3$ and $1<p<2$. Then \eqref{eq:logPS} fails.
\end{proposition}

\begin{proof}
Let $U(y)=c_{0}e^{-|y|^{p/(p-1)}}$, where $c_{0}$ is chosen so that $\int_{\mathbb{R}^N}U^{p}dy=1$. By \cite{PD03, IG}, $U$ gives equality in \eqref{eq:logE}, that is,
\[
\int_{\mathbb{R}^N}U^p\ln U\,dy=\frac{N}{p^{2}}\ln\left(\mathcal{L}_{N,p}A_{E}\right),\qquad A_{E}=\int_{\mathbb{R}^N}|\nabla U|^{p}dy.
\]
Let $u_{\varepsilon}(\rho)=\varepsilon^{-N/p}U(\rho/\varepsilon)\varphi(\rho)$ with $\varphi$ as in the proof of Proposition \ref{prop:fail}, and let $\tilde u_{\varepsilon}=u_{\varepsilon}/M_{\varepsilon}^{1/p}$, where $M_{\varepsilon}=\int_{\mathbb{H}^N}u_{\varepsilon}^{p}dV_{\mathbb{H}}$. Since $U$ and all its derivatives decay faster than any power, the cutoff gives $O(\varepsilon^{k})$ for every $k$, and the term $(\sinh x/x)^{N-1}=1+O(x^{2})$ gives $O(\varepsilon^{2})$ in every integral. So $M_{\varepsilon}=1+O(\varepsilon^{2})$, and since the entropy of $u_\varepsilon$ contains the term $-\frac Np\ln\varepsilon$, the normalization by $M_\varepsilon$ changes it by $O(\varepsilon^{2}|\ln\varepsilon|)$. Hence
\[
\int_{\mathbb{H}^N}\tilde u_{\varepsilon}^p\ln\tilde u_{\varepsilon}\,dV_{\mathbb{H}}=\int_{\mathbb{R}^N}U^p\ln U\,dy-\frac Np\ln\varepsilon+O(\varepsilon^{2}|\ln\varepsilon|),
\qquad
\int_{\mathbb{H}^N}|\nabla_{\mathbb{H}}\tilde u_{\varepsilon}|^{p}dV_{\mathbb{H}}=\varepsilon^{-p}A_{E}\left(1+O(\varepsilon^{2})\right),
\]
and the right side of \eqref{eq:logPS} for $\tilde u_{\varepsilon}$ is equal to
\[
\begin{aligned}
&\frac{N}{p^{2}}\ln(\mathcal{L}_{N,p}A_{E})-\frac Np\ln\varepsilon+\frac{N}{p^{2}}\ln\left(1+O(\varepsilon^{2})-\frac{\beta^{p}}{A_{E}}\varepsilon^{p}\right)\\
&=\int_{\mathbb{H}^N}\tilde u_{\varepsilon}^p\ln\tilde u_{\varepsilon}\,dV_{\mathbb{H}}-\frac{N\beta^{p}}{p^{2}A_{E}}\,\varepsilon^{p}+O(\varepsilon^{2}|\ln\varepsilon|) .
\end{aligned}
\]
Since $p<2$, the right side of \eqref{eq:logPS} is smaller than the left side for small $\varepsilon$.
\end{proof}

Below $p=2$, there is a weaker way to combine the two sharp constants. For $p=1$, Beckner \cite[Theorem 10]{Be15} proved that for $f\in C_{c}^{\infty}(\mathbb{H}^N)$ and $q=\frac{N}{N-1}$,
\begin{equation}\label{eq:beckner}
\left[(N-1)^{N}\|f\|_{L^{1}(\mathbb{H}^N)}^{N}+N^{N}|B|\,\|f\|_{L^{q}(\mathbb{H}^N)}^{N}\right]^{1/N}\leq\int_{\mathbb{H}^N}|\nabla_{\mathbb{H}}f|\,dV_{\mathbb{H}} ,
\end{equation}
and this inequality is sharp. The $N$-th root is exactly what large balls allow, since for large balls the left side is $(N-1)|B_{R}|+O(1)$. For every $N\geq2$ and $1<p<N$, Nguyen \cite[Theorem 4.1]{H18} proved the inequality
\begin{equation}\label{eq:ellN}
\left[\beta^{N}\left(\int_{\mathbb{H}^N}|u|^{p}dV_{\mathbb{H}}\right)^{N/p}+S(N,p)^{N}\left(\int_{\mathbb{H}^N}|u|^{\frac{Np}{N-p}}dV_{\mathbb{H}}\right)^{\frac{N-p}{p}}\right]^{p/N}
\leq\int_{\mathbb{H}^N}|\nabla_{\mathbb{H}}u|^{p}dV_{\mathbb{H}} ,
\end{equation}
which extends the two-dimensional inequalities of Mugelli and Talenti (see \cite{H18}), and which gives \eqref{eq:beckner} as $p\to1$. So the two sharp constants can always be combined with the exponent $\frac Np$. The inequality \eqref{eq:PS} is the case where they are combined with the exponent $1$, and Proposition \ref{prop:fail} shows that this is not possible for $p=1$ and for $1<p<\min\{2,\sqrt N\}$. By Theorem \ref{thm:gap}, it becomes possible when $S(N,p)^{p}$ is replaced by a smaller constant. In fact, the proof of Proposition \ref{prop:fail} gives a little more. For $1<p<\min\{2,\sqrt N\}$ and $s\geq1$, the inequality
\[
\left[\beta^{ps}\left(\int_{\mathbb{H}^N}|u|^{p}dV_{\mathbb{H}}\right)^{s}+S(N,p)^{ps}\left(\int_{\mathbb{H}^N}|u|^{\frac{Np}{N-p}}dV_{\mathbb{H}}\right)^{\frac{s(N-p)}{N}}\right]^{1/s}
\leq\int_{\mathbb{H}^N}|\nabla_{\mathbb{H}}u|^{p}dV_{\mathbb{H}}
\]
fails when $s<\frac{\min\{2,\gamma\}}{p}$, where $\gamma=\frac{N-p}{p-1}$. Indeed, for the functions $u_\varepsilon$ in that proof, the left side minus the right side is at least $c\,\varepsilon^{ps}-C\varepsilon^{\min\{2,\gamma\}}|\ln\varepsilon|$ with $c>0$, which is positive for small $\varepsilon$ when $ps<\min\{2,\gamma\}$.

\medskip
\noindent\textbf{Acknowledgements.} D.~Ganguly was partially supported by the ANRF MATRICS Research Grant (MTR/2023/000331) and the ANRF Advanced Research Grant (ANRF/ARG/2025/000572/MS). N.~Lam was partially supported by an NSERC Discovery Grant. G.~Lu was partially supported by a grant from the Simons Foundation.

\medskip

\noindent{\bf AI Declaration.} The mathematical content of this paper, including all results and proofs, is the work of the authors, who take full responsibility for the content of the paper. AI tools were used only to improve the language of the manuscript after the work was completed.

\end{document}